\documentclass{arxivpaper}

\usepackage{graphicx}
\usepackage{multirow}
\usepackage{commands}
\usepackage{csquotes}
\usepackage{hyperref}
\usepackage{aliascnt}
\usepackage{placeins}
\usepackage[capitalise,noabbrev]{cleveref}
\bibpunct{(}{)}{;}{a}{,}{,}

\usepackage{amsfonts,amsmath}
\usepackage{epstopdf}
\usepackage{tikz}
\usepackage{mathtools}
\usepackage{enumitem}
\usepackage{comment}
\usepackage{booktabs}
\usepackage{subcaption}

\usepackage{siunitx}
\makeatletter

\let\c@algorithm\relax
\makeatother
\usepackage{algorithm}
\usepackage{algpseudocode}
\algrenewcommand\algorithmicrequire{\textbf{Input:}}
\algrenewcommand\algorithmicensure{\textbf{Output:}}
\algnewcommand{\IfThenElse}[3]{
  \State \algorithmicif\ #1\ \algorithmicthen\ #2\ \algorithmicelse\ #3}

\ifpdf
  \DeclareGraphicsExtensions{.eps,.pdf,.png,.jpg}
\else
  \DeclareGraphicsExtensions{.eps}
\fi

\theoremstyle{definition}
\newaliascnt{hypothesis}{theorem}

\aliascntresetthe{hypothesis}
\newaliascnt{assumption}{theorem}
\newtheorem{assumption}[assumption]{\sc Assumption}
\aliascntresetthe{assumption}
\newaliascnt{question}{theorem}
\newtheorem{question}[question]{\sc Question}
\aliascntresetthe{question}
\makeatletter

\let\c@notation\relax

\makeatother
\newaliascnt{notation}{theorem}
\newtheorem{notation}[notation]{\sc Notation}
\aliascntresetthe{notation}

\crefname{hypothesis}{Hypothesis}{Hypotheses}
\Crefname{hypothesis}{Hypothesis}{Hypotheses}
\crefname{assumption}{Assumption}{Assumptions}
\Crefname{assumption}{Assumption}{Assumptions}
\crefname{question}{Question}{Questions}
\Crefname{question}{Question}{Questions}
\crefname{algorithm}{Algorithm}{Algorithms}
\Crefname{algorithm}{Algorithm}{Algorithms}
\crefname{claim}{Claim}{Claims}
\Crefname{claim}{Claim}{Claims}
\crefname{notation}{Notation}{Notations}
\Crefname{notation}{Notation}{Notations}

\usepackage{amsopn}

\DeclarePairedDelimiter{\p}{\lparen}{\rparen}
\DeclarePairedDelimiter{\set}{\{}{\}}
\DeclarePairedDelimiter{\abs}{\lvert}{\rvert}
\DeclarePairedDelimiter{\norm}{\lVert}{\rVert}

\DeclarePairedDelimiterX{\inp}[2]{\langle}{\rangle}{#1, #2}

\begin{document}

\title{Trustworthy AI in numerics: On verification algorithms for neural network-based PDE solvers}
\shorttitle{Verification of neural network-based PDE solvers}

\author{%
{\sc Emil Haugen\thanks{emilhau@math.uio.no}}\\[2pt]
Department of Mathematics, University of Oslo, 0851 Oslo, Norway\\[6pt]
{\sc and}\\[6pt]
Alexei Stepanenko\thanks{as3304@cam.ac.uk}
and Anders C. Hansen\thanks{ach70@cam.ac.uk}\\[2pt]
Department of Applied Mathematics and Theoretical Physics, University of Cambridge, Cambridge CB3 0WA, United Kingdom
}
\shortauthorlist{E. Haugen \emph{et al.}}

\maketitle

\begin{abstract}{%
We present new algorithms for a posteriori verification of neural networks (NNs) approximating solutions to PDEs. We use numerical quadrature to compute upper bounds for $L^2$ norms of NNs and their derivatives. When combined with energy estimates for specific PDEs, this yields verification algorithms which only output approximations with $\varepsilon$-accuracy (in a suitable norm) with respect to the true but unknown solution of the PDE -- for any given $\varepsilon > 0$. This framework enables trustworthy algorithms for NN-based PDE solvers, regardless of training method. Such a posteriori verification is essential because a priori error bounds generally cannot guarantee the accuracy of computed solutions due to the algorithmic undecidability of the optimisation problems used to train NNs.%
}{Trustworthy AI; PDEs; verification algorithms; physics-informed neural networks.}
\end{abstract}

{\small\noindent\textit{2020 Mathematics Subject Classification:}
65M15, 65G20, 65D32, 68T07, 26D10\par}

\section{Introduction}
Fuelled by success in such diverse fields as image analysis, medical imaging, language processing, drug discovery, and computer vision, artificial intelligence (AI) techniques such as deep learning (DL) and neural networks (NNs) are becoming increasingly prevalent in scientific computing. In particular, NNs have emerged as tools for solving partial differential equations (PDEs) \citep{NatureMachIntlEditor_2025} through methods such as physics-informed NNs (PINNs) \citep{dissanayake1994neural-network-based,
      lagaris1998artificial, lagaris2000neural-network,
      raissi2018hidden, raissi2019physics-informed},
 and neural operators \citep{
  Townsend_2024,
 Kovachki_Lanthaler_Mishra_2021,
 Kovachki_Li_Liu_Azizzadenesheli_Bhattacharya_Stuart_Anandkumar_2023,
   Kovachki_Lanthaler_Stuart_2024,
Li_Kovachki_Azizzadenesheli_Liu_Bhattacharya_Stuart_Anandkumar_2021,
 Perdikaris_2021,
 PINO_2024,
 Zha_Biophysics_MICCAI2024,
 Lu_Jin_Pang_Zhang_Karniadakis_2021}.
 In addition to their use in solving differential equations
\citep{Gadapt_2025,
hybridNewton_2025,
deryck2024wpinns,
Grossmann_Komorowska_Latz_Schonlieb_2024},
AI methods have even been proposed in connection with computer-assisted proofs for major open problems in PDEs \citep{wang2023asymptotic, Javier_Discovery_2025}. The most prominent example is of course OpenAI's recently announced resolution of the Navier--Stokes existence and smoothness problem \citep{OpenAI_Navier_Stokes_2026}. Almost simultaneously, PINNs were used to compute a candidate self-similar blowup profile for the Euler equations
\citep{Ganeshram_Duruisseaux_Anandkumar_2026}.

A key challenge in AI for scientific computing and discovery is to guarantee trustworthy results. In fields such as inverse problems and imaging, microscopy, large language models (LLMs), image classification, and computer vision, AI methods can become unstable or provide incorrect results while presenting them as correct (so-called hallucinations)
\citep{bel19, Fawzi2018, fastmri20, TKpaper, Hallucinate_IEEESpec}.
Therefore, the general problem of making AI methods trustworthy and reliable in scientific computing is a major current concern \citep{kutyniok2024mathematics, burger2024learning}. We will address this issue by considering the following basic question.

\begin{question}\label{question:trust-PDE}
    How can AI methods be made trustworthy for PDEs?
\end{question}

By \textit{trustworthy}, we mean that, given some $\eps > 0$, we can guarantee that the approximation produced by the AI method is no more than $\eps$ away from the true solution in some predefined metric.
The primary challenge concerning the trustworthiness of DL methods in solving PDEs lies in the near impossibility of establishing a priori error bounds for trained NNs in practice. This difficulty arises because the training process typically relies on iterative optimisation algorithms designed to minimise a problem-specific loss function.
Rigorously analysing this optimisation procedure -- which involves navigating a complex landscape of hyperparameter choices while tackling highly non-convex objectives -- presents profound theoretical hurdles. Not only does this challenge our analytical capabilities, but it is also -- as recently established -- \emph{algorithmically undecidable} in many cases \citep{Alex1, opt_big, PNAS_SMALE18, gazdagFocm, CRP, AIM, boche2022limitations}.

\subsection{Necessity of a posteriori error bounds -- Undecidability in optimisation}
Recent theoretical results identify scenarios where optimal NNs can be proved to exist, yet no algorithm is capable of computing them reliably \citep{Alex1, MathsOfAdvAttacs, Matt_SIAM_News, gazdagFocm, PNAS_SMALE18, kutyniok2024mathematics, boche2022limitations}. This paradox arises from the following phenomenon. Let $\mathcal{L}: \mathbb{R}^d \rightarrow \mathbb{R}$ be convex and consider the computational problem of finding
\begin{equation}\label{eq:opt}
 \tilde \theta \in \mathrm{argmin}_{\theta \in \mathcal{C}} \mathcal{L}(\theta),
\end{equation}
where $\mathcal{C} \subset \mathbb{R}^d$ is some convex set. For basic objective functions $\mathcal{L}$ and constraints $\mathcal{C}$ such as linear programming, LASSO, and basis pursuit, there are phase transitions in computability/decidability \citep{Alex1, opt_big, PNAS_SMALE18, gazdagFocm, CRP}; see in particular Problem~5 by J.~Lagarias in \citet{AIM}. This occurs when the input is given approximately -- yet with arbitrary precision. Hence, it happens in all cases where floating-point arithmetic, which is the basis of scientific computing, is used.
More precisely, the phenomenon can be described as follows: for an approximation threshold $\eps_1 > 0$, an $\eps$-approximation to a minimiser of \eqref{eq:opt} can be computed for $\eps > \eps_1$, but becomes undecidable (no algorithm exists to compute it) for $\eps < \eps_1$.
This undecidability/non-computability phenomenon will also occur in the non-convex case when computing NNs \citep{MathsOfAdvAttacs}.

\begin{remark}[Non-computability even in the randomised case] Note that the undecidability or non-computability results above \citep{opt_big, PNAS_SMALE18, gazdagFocm, CRP, AIM, BCH_2_2023} also hold in the randomised case. In particular, when $\eps < \eps_1$, there is no randomised algorithm that can compute an $\eps$-approximation with probability higher than $1/2$. This is crucial, as most optimisation algorithms used in modern AI are randomised. Moreover, the approximation threshold $\eps_1$ can be anywhere in the interval $(0,1)$; that is, for any $\eps_1 \in (0,1)$, there is a problem class (of LASSO problems, say) that has a phase transition at $\eps_1$.
\end{remark}

Thus, a priori error bounds for NNs applied to PDEs -- whether stated as general existence results or as characterisations of the minimisers of optimisation problems -- do not, in themselves, guarantee accuracy of the computed network. Such guarantees require bounds that incorporate a detailed analysis of the entire optimisation procedure. However, this level of analysis is only feasible in very special cases, since the problem is, in general, non-computable. This leaves one natural option: a posteriori error bounds. The a posteriori approach is well-established in the classical PDE literature \citep{Verfurth_2013, Carstensen_Merdon_2010, Chamoin_Legoll_2023, Bonito_Canuto_Nochetto_Veeser_2024}, but a general framework for its effective application to DL-based solutions is currently lacking.

Suppose now that we are interested in measuring the error of a given approximation $v$ (the NN) to the true solution $u$ of a PDE in some normed space $X$ with norm $\norm{\cdot}_X$. Ideally, we would like to be able to answer the following question.
\begin{question}\label{question:X-error}
    Given a tolerance $\eps > 0$ and an approximation $v$ of the true but unknown solution $u$, is $\norm{u-v}_X < \eps$?
\end{question}
However, in general, rigorously estimating the desired norm $\norm{u-v}_X$ directly will not be feasible. To resolve this, our residual-based approach relies on constructing inequalities of the form
\begin{equation}\label{eq:XYZ-residuals}
    \norm{u-v}_X \leq C_1 \, \norm{\phi[v]}_Y \leq  C_2 \, \norm{u-v}_Z, 
\end{equation}
where $Y$ and $Z$ are (semi-)normed spaces
and $\phi \colon V \to Y$ is some given operator arising from the PDE with $\phi[u] = 0$, and $V \subset X$ is a subspace where the approximation is computed. Instead of the generally infeasible \cref{question:X-error}, we then consider the following weaker question.

\begin{question}\label{question-sequence}
Can one construct an algorithm with the following property? Suppose $\set{v_m}$ is a sequence of NNs containing a subsequence converging to the unknown true solution $u$ of a given PDE in the normed space $Z$. Given $\eps > 0$, the algorithm determines an approximation $v_m$ such that $\norm{u-v_m}_X < \eps$.
\end{question}

As we show in this paper, \cref{question-sequence} can be answered affirmatively in many cases by applying the framework \eqref{eq:XYZ-residuals}.

 \begin{remark}
In practice, one is not given a sequence as in \cref{question-sequence}, but this formalises the idea that iterative training procedures yield better and better approximations. This is similar to classical approaches such as finite element or finite difference methods, where the convergence analysis typically requires some parameter $h$ to tend to zero, which of course does not happen in practice.

Our verification algorithms will show that for an approximation $v$ we can guarantee that $\norm{u-v}_X < \eps_1$, which is a true but often pessimistic estimate. Indeed, typically $\norm{u-v}_X < \eps_2$ and $\eps_2 \leq \eps_1$, but we will never know $\eps_2$.
\end{remark}

\subsection{Novelty of contributions}
The paper develops new verification algorithms providing rigorous a posteriori error bounds for NNs approximating solutions to PDEs. Our approach is based on residual estimates and rigorous numerical computation of $L^2$ bounds of NNs and their derivatives. PINNs are used for the numerical examples, but our methods work identically for a NN produced by any other training process.  The paper is organised as follows.    
\begin{enumerate}[label={}]
\item \textbf{\cref{sec:NN}}: We introduce the NN model and the residual frameworks for the heat and wave equations. We then construct quadrature rules that we use to rigorously estimate the residual norms. We also formulate the assumptions required on the data and model for the computed error bounds to converge.

\item \textbf{\cref{sec:activation}}:
This section contains the main technical contributions: algorithms producing rigorous bounds for the local variation of derivatives of smooth neural networks with an arbitrary number of layers.

\item \textbf{\cref{sec:applications}}: We derive energy estimates for the heat and wave equations with explicit constants. We then detail how these estimates combine with the algorithms from \cref{sec:activation} to produce rigorous upper bounds converging to the residual norms appearing in the energy estimates. Finally, we use these upper bounds to give informal descriptions of local and global \emph{verification algorithms} that provide an affirmative answer to \cref{question-sequence}.

\item \textbf{\cref{sec:numerical-examples}}: We present numerical examples of the verification algorithms applied to NN-based approximations for the heat and wave equations. The results show that the algorithms provide rigorous a posteriori error bounds for the trained NNs, thus verifying their accuracy with respect to the true solution of the PDEs.

\item \textbf{A programme on trustworthy AI -- Undecidability and computer-assisted proofs included:}
Residual bounds for PDEs often follow from well-posedness results \citep{Keel_Tao_1998, Blair_Smith_Sogge_2009, Bez_Kinoshita_Shiraki_2023}. Combining such theories with rigorous numerical computation of Sobolev norms, as done in this paper, suggests a broader verification programme. However, caution is warranted. Just as undecidability plays a role in optimisation and the construction of NNs, undecidability of the location of solutions to PDEs (see Tao's programme on embedding universal Turing machines into flows of differential equations \citep{Tao_2017, Tao_2016, Tao_averaged_2016}) must play a role in such a theory. Hence, a programme on trustworthy AI in PDEs will touch on numerical analysis, logic and recursion theory, and computer-assisted proofs \citep{Javier_Discovery_2025, Chen_Hou_2025, Chen_Hou_2023, Chen_Hou_2024, Cepelewicz_2022, Wood_2026}.
\end{enumerate}

\subsection{Existing literature}
NN solutions of PDEs, especially through PINNs, have generated a vast literature in recent years; see \citet{cuomo2022scientific} for an overview. PINNs form an \textit{unsupervised learning} framework in which the parameters of an NN solution are chosen in an attempt to minimise discretised residuals given by the PDE.
The literature is heavily focused on applications, but recently progress has also been made on the theoretical side. Studies on the convergence properties and approximation power of NNs in solving PDEs demonstrate the \textit{existence} of accurate approximations
\citep{lu2021priori, yeonjongshin2020convergence, shin2023error,siegel2023greedy, mishra2022estimates,mishra2023estimates, deryck2021approximation}.
Associated probabilistic or asymptotic error estimates are constructed in \citep{deryck2022error, deryck2024error, Zeinhofer_Masri_Mardal_2024}. Distinct from our residual-based approach, \citet{hillebrecht2023rigorous} employ semigroup theory to obtain a posteriori bounds subject to certain assumptions in the steps requiring numerical integration (indeed, our methods for rigorous numerical integration of NNs and their derivatives can be of use here).

In contrast to the comprehensively LLM-based approach to proof construction in \citet{OpenAI_Navier_Stokes_2026}, our work concerns the verification of NNs representing numerical approximations to PDE solutions. Thus it is more closely related to the PINN-based work of \citet{Ganeshram_Duruisseaux_Anandkumar_2026} on blowup for the Euler equations. There, the authors use spline approximations of the trained PINN to obtain interval bounds on derivatives and PDE residuals. The methods developed in \cref{sec:activation} can potentially simplify this certification stage by bounding NN derivatives directly.

\section{Residual verification framework}\label{sec:NN}

This section presents the PDEs which we will consider, and fixes the language and notation we will use to verify their approximate solutions.  

\subsection{Neural-network approximations}
Fix an input dimension $n \in \mathbb{N}$ and let $\cF_{L,w}$ denote scalar neural networks with $L \in \N$ hidden layers, common width $w \in \N$, and smooth activation $\mu=\tanh$. An element $f \in \cF_{L, w}$ acts on  $x\in\R^n$ by
\begin{equation}\label{eqdef:NN}
\begin{aligned}
 z^1(x)&=W^1x+b^1,            & a^1(x)&=\mu(z^1(x)),\\
 z^k(x)&=W^ka^{k-1}(x)+b^k,   & a^k(x)&=\mu(z^k(x)),\quad k=2,\dots,L,\\
 f(x)&=W^{L+1}a^L(x)+b^{L+1}.
\end{aligned}
\end{equation}

\begin{assumption}[Ansatz form]\label{assumption:NN-bounds}
For a bounded spatial domain $\Omega\subset\R^d$ and time interval $[0,T]$, the approximation is given by
\begin{equation}\label{eq:boundary-ansatz}
    v(x,t)=B(x)f(x,t),
\end{equation}
where $f$ is a neural network as in \eqref{eqdef:NN} with $n=d+1$, and $B$ is an explicit smooth function vanishing on $\partial\Omega$.
\end{assumption}
In our applications with $\Omega=(0,1)^d$, the standard product factor enforcing homogeneous Dirichlet data is $B(x)=\prod_{i=1}^d x_i(1-x_i)$.

\subsection{The residual verification principle}\label{subsec:residual-principle}
Let $u$ be the unique, unknown solution to a given PDE and let $v$ be the proposed approximation. A residual estimate has the schematic form
\begin{equation}\label{eq:posteriori-residuals}
    \norm{u-v}_{X} \leq C\,\mathcal R(v;\text{data}),
\end{equation}
where $\cR$ is a continuous functional. For linear parabolic and hyperbolic problems, $\cR$ is a sum of $L^2$ norms. If the constant $C$ is known and $\mathcal R$ can be rigorously bounded, then the unknown error has a computable upper bound.

We will consider the source-free Dirichlet heat and wave equations as model problems:
\begin{equation}\label{eq:heat-wave-model-framework}
\begin{array}{@{}c@{\qquad}c@{}}
\text{Heat} & \text{Wave} \\[0.3em]
\begin{cases}
\partial_t u-\kappa\Delta u=0 & \text{in }\Omega\times(0,T],\\
u=0 & \text{on }\partial\Omega\times[0,T],\\
u=g & \text{on }\Omega\times\set{t=0},
\end{cases}
&
\begin{cases}
\partial_t^2 u-c^2\Delta u=0 & \text{in }\Omega\times(0,T],\\
u=0 & \text{on }\partial\Omega\times[0,T],\\
u=g,\quad \partial_tu=h & \text{on }\Omega\times\set{t=0}.
\end{cases}
\end{array}
\end{equation}
For the heat and wave equations, we have the respective estimates establishing the framework \eqref{eq:posteriori-residuals}:
\begin{align}
    \norm{u-v}
    &\leq
    C\left(
    \norm{g-v(\cdot,0)}
    +\norm{(\partial_t-\kappa\Delta)v}_{L^2(\Omega\times(0,T))}
    \right),
    \label{eq:heat-h1-framework}\\
    \norm{u-v}
    &\leq
    C\left(
    \norm{g-v(\cdot,0)}
    +\norm{h-\partial_t v(\cdot,0)}
    +\norm{(\partial_t^2-c^2\Delta)v}_{L^2(\Omega\times(0,T))}
    \right).
    \label{eq:wave-h1-framework}
\end{align}
We refer to the norms involving $g, h$ as \emph{(initial) data residuals} and the norms involving the equation operator as \emph{PDE residuals}. The value of $C$ depends on the equation and domain, and the precise norms used for the left-hand side and the data residuals depend on the regularity of the initial data. For example, for the heat equation, if we only assume $g\in L^2(\Omega)$ and use that norm on the RHS of \eqref{eq:heat-h1-framework}, we get an $L^\infty(0,T;L^2(\Omega))$ norm on the LHS. If, however, we assume $g\in H^1_0(\Omega)$ and use the $H^1$ norm on $g - v(\cdot,0)$, we get a much stronger norm on the LHS and hence a stronger error bound. We will give precise statements and proofs with explicit constants in \cref{sec:applications}.

\begin{remark}[Removing the hard boundary ansatz]
The product form $v=Bf$ is a convenient way to enforce the homogeneous Dirichlet
condition exactly, and it is the ansatz used in the numerical experiments below.
It is not intrinsic to the residual-verification framework.  Let $A$ denote the PDE operator.  One can verify
a pure NN approximation $v=f$ by introducing a function $b$ with $b=f$ on
$\partial\Omega\times[0,T]$ and applying the corresponding homogeneous
Dirichlet estimate to $f-b$.  The resulting estimate then contains the norm
$\norm{b}_{X}$, the modified initial residuals, and the
PDE residual $A(f-b)$.

Alternatively, $b$ can be chosen to solve the
auxiliary homogeneous problem $Ab=0$, with boundary data
$f|_{\partial\Omega\times[0,T]}$ and compatible auxiliary initial data when
required.  For the heat equation, this means choosing $b(\cdot,0)=b_0$ with
$b_0|_{\partial\Omega}=f(\cdot,0)|_{\partial\Omega}$. Choosing
$Ab=0$ removes the contribution from the PDE residual, but $b$ is then
defined as the solution of an auxiliary PDE, so its values and norms cannot
usually be evaluated by the same direct grid-point formulas used for explicit
NN residuals, and must instead be separately bounded and verified.
\end{remark}

\subsection{Quadrature setup and assumptions}\label{subsec:quadrature}
Given the energy estimates \eqref{eq:heat-h1-framework} and \eqref{eq:wave-h1-framework}, it is clear that to bound the approximation error of $v$ without knowing $u$, we must compute rigorous upper bounds on norms involving the initial data and the approximant $v$ along with its derivatives. We now detail which assumptions are made on the domain and data, and what we will require of the neural network-based approximants, in order to construct the desired norm bounds.

\subsubsection{Domain geometry and partitioning}

\begin{assumption}[Rectangular domain geometry]\label{assumption:rectangular-domain}
  Let $\Omega \subset \R^d$ be a bounded domain, let $T>0$, and write $\Omega_T \coloneqq \Omega \times (0, T]$. Fix $\boldsymbol{\eps}_x\in(0,\infty)^d$ and
  $\eps_t>0$. Set
  $\boldsymbol{\eps}\coloneqq(\boldsymbol{\eps}_x,\eps_t)\in(0,\infty)^{d+1}$.
  For a radius vector
  $\boldsymbol r=(\boldsymbol r_x,r_t)\in[0,\infty)^{d+1}$,
  $y\in\R^d$, and $t\in\R$, define
  \begin{equation}\label{quad-spacetime-box}
    \begin{aligned}
      \mathcal B_{\boldsymbol r_x}(y)
      &\coloneqq
      \set{x\in\R^d: |x_j-y_j|\leq r_{x,j},\ j=1,\dots,d},\\
      I_{r_t}(t)&\coloneqq[t-r_t,t+r_t],\\
      \mathcal B_{\boldsymbol r}(y,t)
      &\coloneqq
      \mathcal B_{\boldsymbol r_x}(y)\times I_{r_t}(t).
    \end{aligned}
  \end{equation}
  We assume that $\Omega$ is rectangular in the sense that there exist $y_k\in\Omega$ for $k=1,\dots,N_{\mathrm{space}}$ and
  $t_m\in(0, T)$ for $m=1,\dots,N_{\mathrm{time}}$ such that the spatial boxes $\mathcal B_{\boldsymbol{\eps}_x}(y_k)$ and the time intervals $I_{\eps_t}(t_m)$ are pairwise internally disjoint and satisfy
  \[
    \overline{\Omega}
    =
    \bigcup_{k=1}^{N_{\mathrm{space}}}\mathcal B_{\boldsymbol{\eps}_x}(y_k),
    \qquad
    [0,T]
    =
    \bigcup_{m=1}^{N_{\mathrm{time}}} I_{\eps_t}(t_m).
  \]
  Consequently, the rectangles $\mathcal B_{\boldsymbol{\eps}}(y_k,t_m)$ are pairwise internally disjoint and
  \[
    \overline{\Omega_T}
    =
    \overline{\Omega}\times[0,T]
    =
    \bigcup_{k, m}
    \mathcal B_{\boldsymbol{\eps}}(y_k,t_m).
  \]
  We write
  \(
    \mathcal P_\Omega
    \coloneqq
    \set*{\mathcal B_{\boldsymbol{\eps}_x}(y_k)}_k,
  \) and
  \(  \mathcal P
    \coloneqq
    \set*{\mathcal B_{\boldsymbol{\eps}}(y_k,t_m)}_{k,m}.
  \)
\end{assumption}

This assumption is not as restrictive as it may appear, since we are interested in \emph{upper bounds} for norms, and any bounded domain can be covered by rectangles. If the residual functions can be extended smoothly (e.g. by reflections of the data) to a rectangular cover, then norm bounds computed on the cover will be valid upper bounds, and convergence can be ensured by letting the cover approach the domain. For ease of notation, we assume a common radius vector $\boldsymbol{\eps}$ for all cells, but it is straightforward to allow for different radii; indeed, such adaptive partitions are implemented in our code.

\subsubsection{Quadrature rules and error bounds}
\label{subsubsec:cellwise-quadrature-bounds}
We now introduce the quadrature methods and notation that we will use to bound the residual norms of \eqref{eq:heat-h1-framework} and \eqref{eq:wave-h1-framework} in \cref{sec:applications}. We focus on $L^2$ norms since they are central to our applications, but it is straightforward to extend our work to general $L^p$ norms.

Given \cref{assumption:rectangular-domain}, let
\[
  C=\zeta_C+\prod_{q=1}^s[-\eps_{C,q},\eps_{C,q}]
\]
be a cell centred at $\zeta_C \in \R^s$ with radius vector $\boldsymbol{\eps}_C\in(0,\infty)^s$, where $s=d$ for spatial cells and $s=d+1$ for spacetime cells. For a scalar integrand $\phi \colon \R^s \to \R$, we consider the local approximations
\begin{equation}\label{eq:cell-polynomial-models}
  P_C^0(\zeta)\coloneqq \phi(\zeta_C),
  \qquad
  P_C^1(\zeta)\coloneqq
  \phi(\zeta_C)+\nabla\phi(\zeta_C)\cdot(\zeta-\zeta_C),
  \qquad \zeta\in C.
\end{equation}
Define the local quadrature estimates
\begin{equation}\label{eq:cell-quadrature-rule}
  Q_C^r(\phi)
  \coloneqq
  \int_C P_C^r(\zeta)^2\,d\zeta, \quad r\in \{0,1\}.
\end{equation}
These can be computed explicitly from the cell centre and radius vector,
\begin{equation}\label{eq:cell-quadrature-values}
  Q_C^0(\phi)=|C|\,\phi(\zeta_C)^2,
  \qquad
  Q_C^1(\phi)
  =
  |C|\left(\phi(\zeta_C)^2
  +\frac{1}{3}\sum_{q=1}^s\eps_{C,q}^2\partial_q\phi(\zeta_C)^2\right).
\end{equation}
For sufficiently smooth integrands, both of these rules can be seen to have quadratic convergence rates with respect to the cell diameter. However, we will shortly see that for the actual error bounds computed \emph{in practice}, the classical midpoint rule $Q^0$ under $C^1$ regularity assumptions gives linear convergence, while the affine rule $Q^1$ remains quadratic, at the cost of requiring $C^2$ regularity.

For any rectangular spatial or spacetime partition $\mathcal T$, define the global quadrature estimate
\begin{equation}\label{eq:partition-quadrature-rule}
  Q_{\mathcal T}^r(\phi)\coloneqq \sum_{C\in\mathcal T}Q_C^r(\phi),
  \qquad r\in\{0,1\}.
\end{equation}
Suppose for the moment that we have computable bounds on the first and second derivatives:
\begin{equation}\label{eq:cell-derivative-coefficients}
  b_{q,C}(\phi)\geq\sup_{\zeta\in C}|\partial_q\phi(\zeta)|,
  \qquad
  h_{q\ell,C}(\phi)\geq
  \sup_{\zeta\in C}|\partial_{q\ell}\phi(\zeta)|.
\end{equation}
Writing $\rho_q=|\zeta_q-\zeta_{C,q}|$, define
\begin{equation}\label{eq:cell-remainder-bounds}
  D_C^0(\phi)(\zeta)
  \coloneqq\sum_{q=1}^s b_{q,C}(\phi)\rho_q,
  \qquad
  D_C^1(\phi)(\zeta)
  \coloneqq\frac12\sum_{q,\ell=1}^s
  h_{q\ell,C}(\phi)\rho_q\rho_\ell.
\end{equation}
By Taylor's theorem, 
\begin{equation}\label{eq:cell-model-errors}
  |\phi(\zeta)-P_C^r(\zeta)|\leq D_C^r(\phi)(\zeta),
  \qquad \zeta\in C,\quad r\in\{0,1\}.
\end{equation}
We also set
\begin{equation}\label{eq:cell-polynomial-bounds}
  A_C^0(\phi)(\zeta)\coloneqq|\phi(\zeta_C)|,
  \qquad
  A_C^1(\phi)(\zeta)\coloneqq
  |\phi(\zeta_C)|+\sum_{q=1}^s
  |\partial_q\phi(\zeta_C)|\rho_q,
\end{equation}
so that $|P_C^r|\leq A_C^r(\phi)$. Writing $\phi^2-(P_C^r)^2=2P_C^r(\phi-P_C^r)+(\phi-P_C^r)^2$ gives the local bound
\begin{equation}\label{eq:cell-verified-quadrature-error}
  \left|
  \int_C \phi(\zeta)^2\,d\zeta-Q_C^r(\phi)
  \right|
  \leq
  2\int_C A_C^r(\phi)D_C^r(\phi)\,d\zeta
  +\norm{D_C^r(\phi)}_{L^2(C)}^2
  \eqqcolon E_C^r(\phi).
\end{equation}
Then $E_C^r(\phi)$ can be evaluated analytically using the integrals
\begin{equation}\label{eq:centered-cell-moments}
  \mathfrak M_\alpha(C)
  \coloneqq\int_C\prod_{q=1}^s
  |\zeta_q-\zeta_{C,q}|^{\alpha_q}\,d\zeta
  =\prod_{q=1}^s\frac{2\eps_{C,q}^{\alpha_q+1}}{\alpha_q+1},
  \qquad \alpha\in\N_0^s.
\end{equation}
In particular, writing $a_C=|\phi(\zeta_C)|$ and $c_{q,C}=|\partial_q\phi(\zeta_C)|$, we have the computable expressions 
\begin{align}
 E_C^0(\phi)
 &=2a_C\sum_{q=1}^s b_{q,C}(\phi)\mathfrak M_{e_q}(C)
 +\sum_{q,\ell=1}^s b_{q,C}(\phi)b_{\ell,C}(\phi)
 \mathfrak M_{e_q+e_\ell}(C),
 \label{eq:q0-cell-moment-correction}\\
 E_C^1(\phi)
 &=\sum_{q,\ell=1}^s h_{q\ell,C}(\phi)
 \left(a_C\mathfrak M_{e_q+e_\ell}(C)
 +\sum_{i=1}^s c_{i,C}\mathfrak M_{e_q+e_\ell+e_i}(C)\right)
 \notag\\
 &\quad+\frac14\sum_{q,\ell,m,n=1}^s
 h_{q\ell,C}(\phi)h_{mn,C}(\phi)
 \mathfrak M_{e_q+e_\ell+e_m+e_n}(C).
 \label{eq:q1-cell-moment-correction}
\end{align}
On such a partition $\mathcal T$, the global error is bounded by
\begin{equation}\label{eq:partition-verified-quadrature-error}
  \left|
  \int_{\bigcup_{C\in\mathcal T}C}\phi(\zeta)^2\,d\zeta
  - Q_{\mathcal T}^r(\phi)
  \right|
  \leq
  \sum_{C\in\mathcal T}E_C^r(\phi)
  \eqqcolon   E_{\mathcal T}^r(\phi)
  \quad r\in\{0,1\}.
\end{equation}
Consequently,
\begin{equation}\label{eq:global-norm-bound}
  \norm{\phi}_{L^2(\bigcup_{C\in\mathcal T}C)}
  \leq
  \left(Q_{\mathcal T}^r(\phi)+E_{\mathcal T}^r(\phi)\right)^{1/2}
  \eqqcolon\mathcal U_{\mathcal T}^r(\phi).
\end{equation}

We must ensure that the computed upper bounds $\mathcal{U}_{\mathcal T}^r$ actually converge to the norms they estimate. This is a stronger statement than the simple assertion that the quadrature estimates themselves converge.

\begin{lemma}[Convergence of verified norm bounds]
\label{lem:verified-norm-bound-convergence}
Let $G=\Omega$ and $s=d$, or let $G=\Omega_T$ and $s=d+1$.
Let $r\in\{0,1\}$, $\phi\in C^{r+1}(\overline G)$, and suppose 
$\mathcal T$ is a rectangular partition of $G$ as in
\cref{assumption:rectangular-domain}. Set 
\[
 \eps_{\mathcal T}\coloneqq
 \max_{C\in\mathcal T}\max_{1\leq q\leq s}\eps_{C,q}.
\]
Suppose the coefficients in \eqref{eq:cell-derivative-coefficients} used to
construct $\mathcal U_{\mathcal T}^r(\phi)$ satisfy the following bounds for
fixed constants $\overline b_q$ and $\overline h_{q\ell}$ independent of
$\mathcal T$:
\begin{equation}\label{eq:uniform-bound-coeffs}
   b_{q,C}(\phi)\leq\overline b_q\quad\text{if }r=0,
   \qquad
   h_{q\ell,C}(\phi)\leq\overline h_{q\ell}\quad\text{if }r=1,
\end{equation}
for every $C\in\mathcal T$ and all $q, \ell \in \set{1, \dots, s}$. Then there is a
constant $K_\phi>0$, independent of $\mathcal T$, such that
\begin{equation}\label{eq:verified-norm-bound-convergence}
 0\leq
 \mathcal U_{\mathcal T}^r(\phi)-\norm{\phi}_{L^2(G)}
 \leq K_\phi \eps_{\mathcal T}^{r+1}.
\end{equation}
\end{lemma}

\begin{proof}
Set
\[
 \delta_{\mathcal T}^r(\phi)
 \coloneqq
 \left(
 \sum_{C\in\mathcal T}\norm{D_C^r(\phi)}_{L^2(C)}^2
 \right)^{1/2}.
\]
Define the piecewise polynomial $P_{\mathcal T}^r$ by
$P_{\mathcal T}^r|_C=P_C^r$.  The reverse triangle inequality and
\eqref{eq:cell-model-errors} give
\begin{equation}\label{eq:generic-quadrature-comparison}
 \left|
 Q_{\mathcal T}^r(\phi)^{1/2}-\norm{\phi}_{L^2(G)}
 \right|
 =\left|
 \norm{P_{\mathcal T}^r}_{L^2(G)}-\norm{\phi}_{L^2(G)}
 \right|
 \leq\delta_{\mathcal T}^r(\phi).
\end{equation}

For $r=0$, the definitions give
$\norm{A_C^0(\phi)}_{L^2(C)}^2=Q_C^0(\phi)$.  For $r=1$, writing
$a_C=|\phi(\zeta_C)|$, $c_{q,C}=|\partial_q\phi(\zeta_C)|$, and
$\rho_q=|\zeta_q-\zeta_{C,q}|$, the Cauchy--Schwarz inequality in
$\R^{s+1}$ gives
\[
 A_C^1(\phi)(\zeta)^2
 \leq(s+1)\left(a_C^2+\sum_{q=1}^s c_{q,C}^2\rho_q^2\right).
\]
Integration and \eqref{eq:cell-quadrature-values} therefore give
$\norm{A_C^1(\phi)}_{L^2(C)}^2\leq(s+1)Q_C^1(\phi)$.  Hence, for either rule,
\begin{equation}\label{eq:generic-polynomial-bound}
 \left(
 \sum_{C\in\mathcal T}\norm{A_C^r(\phi)}_{L^2(C)}^2
 \right)^{1/2}
 \leq c_{s,r}Q_{\mathcal T}^r(\phi)^{1/2},
 \qquad c_{s,r}\coloneqq\sqrt{1+rs}.
\end{equation}
By the Cauchy--Schwarz inequality in $L^2(G)$,
the definition \eqref{eq:cell-verified-quadrature-error}, and
\eqref{eq:generic-polynomial-bound},
\[
 E_{\mathcal T}^r(\phi)
 \leq
 2c_{s,r}Q_{\mathcal T}^r(\phi)^{1/2}
 \delta_{\mathcal T}^r(\phi)
 +\delta_{\mathcal T}^r(\phi)^2.
\]
Since $c_{s,r}\geq1$, it follows by adding $Q_{\mathcal T}^r(\phi)$ to both sides that
\[
 \mathcal U_{\mathcal T}^r(\phi)
 \leq
 Q_{\mathcal T}^r(\phi)^{1/2}
 +c_{s,r}\delta_{\mathcal T}^r(\phi).
\]
Together with \eqref{eq:global-norm-bound} and
\eqref{eq:generic-quadrature-comparison}, this proves
\begin{equation}\label{eq:generic-convergence-from-above}
 0\leq
 \mathcal U_{\mathcal T}^r(\phi)-\norm{\phi}_{L^2(G)}
 \leq(c_{s,r}+1)\delta_{\mathcal T}^r(\phi).
\end{equation}

Finally, the uniform bounds \eqref{eq:uniform-bound-coeffs} and
\eqref{eq:cell-remainder-bounds} imply
\[
 \delta_{\mathcal T}^r(\phi)
 \leq |G|^{1/2}
 \begin{cases}
 \displaystyle
 \eps_{\mathcal T}\sum_{q=1}^s\overline b_q,&r=0,\\[1ex]
 \displaystyle
 \frac{\eps_{\mathcal T}^2}{2}
 \sum_{q,\ell=1}^s\overline h_{q\ell},&r=1.
 \end{cases}
\]
\end{proof}

\begin{remark}[Vector-valued residuals]\label{rem:vector-residual-bounds}
For a vector-valued residual
$\boldsymbol\phi=(\phi_1,\dots,\phi_p)$, define
\[
 \mathcal U_{\mathcal T}^r(\boldsymbol\phi)
 \coloneqq
 \left(
 \sum_{\ell=1}^p\mathcal U_{\mathcal T}^r(\phi_\ell)^2
 \right)^{1/2}.
\]
Since
\[
  \norm{\boldsymbol\phi}_{L^2}^2
  =\sum_{\ell=1}^p\norm{\phi_\ell}_{L^2}^2,
\]
this is a verified upper bound obtained by applying the scalar construction to
each component.  Moreover, if component $\ell$ satisfies
\eqref{eq:verified-norm-bound-convergence} with constant $K_\ell$, then the
reverse triangle inequality in $\R^p$ gives
\[
 0\leq
 \mathcal U_{\mathcal T}^r(\boldsymbol\phi)
 -\norm{\boldsymbol\phi}_{L^2}
 \leq
 \left(\sum_{\ell=1}^pK_\ell^2\right)^{1/2}
 \eps_{\mathcal T}^{r+1}.
\]
\end{remark}

\subsubsection{Regularity assumptions for quadrature}

We now give the assumptions needed to compute the local bounds in
\eqref{eq:cell-derivative-coefficients}.  The rule $Q^0$ uses first-derivative
bounds for $\phi$, whereas $Q^1$ uses second-derivative bounds.
Since the PDE residuals involve second derivatives of the
approximant, the $Q^1$ quadrature requires derivative bounds
up to order four.  This is encoded in the following assumption for the NN.
\begin{assumption}[Network derivative bounds]\label{ass:network-derivative-bounds}
  For a given neural network $f \in \cF_{L, w}$ and each multi-index $\alpha \in \N_0^{d+1}$ with
  $\abs{\alpha} \leq 4$, we have computable functions
  $\mathcal{E}^{\alpha}_f \colon \R^{d+1} \times [0, \infty)^{d+1} \to \R_+$ such that
  \[
    \abs{\partial^\alpha f(x, s)}
    \leq
    \mathcal{E}^\alpha_f((y, t),\boldsymbol{\eps})
    \quad
    \text{for all }
    \quad (x, s) \in \mathcal B_{\boldsymbol{\eps}}(y,t).
  \]
\end{assumption}

The main technical challenge is replacing \cref{ass:network-derivative-bounds}
with \emph{theorems} providing those bounds, which is what we will do in
\cref{sec:activation} using our novel algorithms.  For the explicit boundary
factor $B$ we make analogous assumptions.

\begin{assumption}[Boundary factor regularity]\label{ass:B-moduli}
For every multi-index $\beta$ with $|\beta|\leq 4$, the point values of
$\partial^\beta B$ are available.  We are also given nondecreasing moduli
$\omega_B^\beta\colon[0,\infty)\to[0,\infty)$, continuous and vanishing at
$0$, such that, for every $x,y\in\overline{\Omega}$,
\begin{equation}\label{eq:B-moduli}
  |\partial^\beta B(x)-\partial^\beta B(y)|
  \leq
  \omega_B^\beta\p*{\norm{x-y}_\infty}.
\end{equation}
\end{assumption}

The regularity available for the prescribed
data determines which residual norm choices and quadrature rules are admissible.

\begin{assumption}[Initial data regularity]\label{ass:data-derivative-bounds}
Write
\begin{equation}\label{eq:initial-data}
  v_0(x)\coloneqq v(x,0),
  \qquad
  \dot v_0(x)\coloneqq \partial_t v(x,0).
\end{equation}
The data residuals are built from $g-v_0$ and, for the wave equation,
$h-\dot v_0$.  For each initial datum $a\in\{g,h\}$ that appears, let $k_a$ be
the largest order of spatial derivative of $a$ occurring in the chosen data
residual.  If that residual is to be computed with quadrature rule $Q^r$,
assume $a\in C^{k_a+r+1}(\overline\Omega)$ and that point values of 
$\partial^\beta a$ are available for every multi-index $\beta$ with
$|\beta|\leq k_a+r$.  For every multi-index $\beta$ with
$|\beta|\leq k_a+r+1$, we are given a computable function
\(
  \mathcal D_a^\beta\colon
  \overline\Omega\times[0,\infty)^d\to\R_+
\)
that is locally bounded and satisfies
\begin{equation}\label{eq:data-derivative-bounds}
  |\partial^\beta a(x)|
  \leq
  \mathcal D_a^\beta(y,\boldsymbol\eps),
  \quad x\in
  \mathcal B_{\boldsymbol\eps}(y)\cap\overline\Omega.
\end{equation}
\end{assumption}

\section{Derivative bounds for neural networks}\label{sec:activation}

This section is devoted to the technical challenge of replacing \cref{ass:network-derivative-bounds} with theorems that provide the computable bound functions $\mathcal{E}^\alpha_f$.

\subsection{Setup: the local derivative bound problem}
Throughout this section, we fix $f\in\cF_{L,w}$ as in \eqref{eqdef:NN}.
For $y\in\R^n$ and a radius vector
$\boldsymbol{\eps}=(\eps_1,\dots,\eps_n)\in[0,\infty)^n$, write
\[
    \mathcal B_{\boldsymbol{\eps}}(y)
    \coloneqq
    \{x\in\R^n:|x_j-y_j|\le \eps_j,\ j=1,\dots,n\}.
\]
We allow $\mathcal B_{\boldsymbol{\eps}}(y)$ to
be degenerate when one or more components of $\boldsymbol{\eps}$ vanish.
In our model problems, we will have $n=d+1$, but we keep the notation general throughout this section.
The goal of this section is to construct, for $m \in \N$, an algorithm
$\alg{m}$ that takes as input a multi-index $\alpha \in \N_0^n$ with
$\abs{\alpha}=m$, a rectangle given by its centre $y\in\R^n$ and
radius vector $\boldsymbol{\eps}\in[0,\infty)^n$, and the NN $f$
(represented by its parameters). The algorithm should compute a number
$\alg{|\alpha|}(\alpha,f,y,\boldsymbol{\eps})$ satisfying
\begin{equation}\label{eq:derivative-difference-bound}
    |\partial^\alpha f(x)-\partial^\alpha f(y)|
    \leq \alg{|\alpha|}(\alpha,f,y,\boldsymbol{\eps}),
    \qquad x\in\mathcal B_{\boldsymbol{\eps}}(y).
\end{equation}
This will achieve the goal of turning \cref{ass:network-derivative-bounds} into a theorem, since then
\begin{equation}\label{eq:derivative-bound}
    |\partial^\alpha f(x)|
    \leq |\partial^\alpha f(y)|+\alg{|\alpha|}(\alpha,f,y,\boldsymbol{\eps})
    \eqqcolon \mathcal{E}^\alpha_f(y,\boldsymbol{\eps}),
    \qquad x\in\mathcal B_{\boldsymbol{\eps}}(y).
\end{equation}

\subsubsection{Common algorithmic structures and notation}\label{sec:common-alg-structure}
Before constructing the specific derivative-bound algorithms, we fix notation and present their common structure. We consistently use the notation of \eqref{eqdef:NN} when referring to the quantities computed throughout the layers of the NN. Since the output layer can be seen as just another affine layer, we write $z^{L+1}\coloneqq f$ when notationally convenient. If $v, u\in\R^w$ and $M\in\R^{w\times n}$, then $v\odot u$ is the Hadamard product, and $v\odot M$ is the Hadamard product of $v$ with every column of $M$. Vector or matrix inequalities, absolute values, and derivatives are componentwise. The vector in $\R^n$ with entries all equal to $1$ is denoted $\mathbf{1}_n$.

For any nonzero multi-index $\alpha$, the idea behind constructing the right-hand side of \eqref{eq:derivative-difference-bound} is to consider the functions
\begin{equation}\label{eqdef:p-h-alpha}
    p_\alpha^k \coloneqq \partial^\alpha z^k,
    \qquad
    h_\alpha^k \coloneqq \partial^\alpha a^k.
\end{equation}
At the output layer, $z^{L+1}=f$, we similarly write
$p_\alpha^{L+1}\coloneqq\partial^\alpha f$.  Also note the relation
\begin{equation}\label{eq:p-h-recurrence}
    p_\alpha^{k+1}
    =
    W^{k+1}h_\alpha^k,
    \qquad k=1,\dots,L.
\end{equation}
The algorithms then construct vectors $P_\alpha^k, B_\alpha^k \in \R_+^w$ satisfying, for $k = 1,\dots,L$ and every $x\in\mathcal B_{\boldsymbol{\eps}}(y)$,
\begin{subequations}\label{eq:variation-bounds}
\begin{align}
    P_\alpha^k &\ge |p_\alpha^k(x)-p_\alpha^k(y)|, \label{eq:preactivation-variation-bound}\\
    B_\alpha^k &\ge |h_\alpha^k(x)-h_\alpha^k(y)|. \label{eq:activation-variation-bound}
\end{align}
\end{subequations}
The final bound \eqref{eq:derivative-difference-bound} follows directly from \eqref{eq:activation-variation-bound} at layer $L$,
\begin{equation}\label{eq:output-variation-bound}
    |p_\alpha^{L+1}(x)-p_\alpha^{L+1}(y)|
    =
    \abs*{\partial^\alpha f(x)-\partial^\alpha f(y)}
    \le
    |W^{L+1}|B_\alpha^L \eqqcolon P^{L+1}_\alpha.
\end{equation}
When $\alpha$ is of order one, $\alpha = e_\ell$, we will also simply use the subscript $\ell$ in the above quantities. To compute the bounds \eqref{eq:activation-variation-bound} for a given $\alpha$, we will make repeated use of the following technical lemma with $2\le s\le |\alpha|+1$.

\begin{lemma}[Product-difference bound]\label{lem:xi-lemma}
Let $s \in \N$.  For each $i=1,\dots,s$, suppose $\xi_i,\xi_{i,0}, E_i\in\R^w$ satisfy $|\xi_i-\xi_{i,0}|\le E_i$.  Then
\begin{equation}\label{eq:xi-lemma}
\left|\bigodot_{i=1}^s \xi_i-\bigodot_{i=1}^s \xi_{i,0}\right|
\le
\bigodot_{i=1}^s(|\xi_{i,0}|+E_i)-\bigodot_{i=1}^s|\xi_{i,0}| .
\end{equation}
\end{lemma}
\begin{proof}
See \cref{app:xi-lemma-proof}.
\end{proof}

We also need a computational formula for the activation derivatives. 
\begin{lemma}[General chain rule for activations]\label{lem:faa-di-bruno-expansion}
Fix a multi-index of size $m\ge1$,
$\alpha=e_{\ell_1}+\cdots+e_{\ell_m}$.  Let
$\Pi_\alpha$ be the set of all partitions of $\{1,\dots,\abs{\alpha}\}$ into nonempty
blocks, and for each nonempty $I\subseteq\{1,\dots,\abs{\alpha}\}$ define
$\beta_I\coloneqq\sum_{i\in I}e_{\ell_i}$.  Then for each hidden layer $k=1,\dots,L$,
\begin{equation}\label{eq:abstract-product-expansion}
    h_\alpha^k
    =
    \sum_{\pi\in\Pi_\alpha}
    \mu^{(|\pi|)}(z^k)
    \odot
    \bigodot_{I\in\pi} p_{\beta_I}^k.
\end{equation}
\end{lemma}
\begin{proof}
This follows from the ordinary Faà di Bruno formula applied componentwise to
$a^k=\mu\circ z^k$.
\end{proof}

The chain rule expansion lets us formulate the general procedure for propagating the bounds \eqref{eq:variation-bounds} through the network.  

\begin{lemma}[Bound propagation via chain rule]\label{lem:bound-iteration}
Let $\alpha$ be a nonzero multi-index, and let $\Pi_\alpha$ and $\beta_I$ be
as in \cref{lem:faa-di-bruno-expansion}.  Set
\[
    \Gamma_\alpha\coloneqq
    \{\beta_I:\pi\in\Pi_\alpha,\ I\in\pi\}.
\]
Suppose that, for a given layer
$\kappa \in \set{1,\dots, L}$, we have
$P_\gamma^\kappa,Q^{|\pi|,\kappa}\in\R_+^w$ satisfying for every
$x\in\mathcal B_{\boldsymbol{\eps}}(y)$,
\begin{subequations}\label{eq:bounds-p-and-mu}
\begin{align}
    \left|p_\gamma^\kappa(x)-p_\gamma^\kappa(y)\right|
    &\le P_\gamma^\kappa,
    &&\gamma\in\Gamma_\alpha,\label{eq:bounds-p-and-mu-a}\\
    \left|\mu^{(|\pi|)}(z^\kappa(x))-\mu^{(|\pi|)}(z^\kappa(y))\right|
    &\le Q^{|\pi|,\kappa},
    &&\pi\in\Pi_\alpha.\label{eq:bounds-p-and-mu-b}
\end{align}
\end{subequations}
Define
\begin{equation}\label{eqdef:phat-Mk}
    \widehat p_\gamma^\kappa\coloneqq |p_\gamma^\kappa(y)|+P_\gamma^\kappa,
    \qquad
    M_j^\kappa\coloneqq |\mu^{(j)}(z^\kappa(y))|+Q^{j,\kappa}
\end{equation}
and set
\begin{equation}\label{eq:abstract-B-alpha}
    B_\alpha^\kappa
    \coloneqq
    \sum_{\pi\in\Pi_\alpha}
    \left(
    M_{|\pi|}^\kappa
    \odot
    \bigodot_{I\in\pi}\widehat p_{\beta_I}^\kappa
    -
    \left|
    \mu^{(|\pi|)}(z^\kappa(y))
    \odot
    \bigodot_{I\in\pi}p_{\beta_I}^\kappa(y)
    \right|
    \right).
\end{equation}
Then at layer $\kappa+1$ we have the following inequality, which lets us propagate the bounds \eqref{eq:variation-bounds} through the network:
\begin{equation}\label{eq:affine-successor-variation-bound}
    |p_\alpha^{\kappa+1}(x)-p_\alpha^{\kappa+1}(y)|
    \le |W^{\kappa+1}|B_\alpha^\kappa,
    \qquad x\in\mathcal B_{\boldsymbol{\eps}}(y).
\end{equation}

\end{lemma}
\begin{proof}
Fix $x\in\mathcal B_{\boldsymbol{\eps}}(y)$ and a partition $\pi\in\Pi_\alpha$.  Write
$\pi=(I_1,\dots,I_{|\pi|})$.  To apply
\cref{lem:xi-lemma}, set
\[
    \xi_1=\mu^{(|\pi|)}(z^\kappa(x)),
    \qquad
    \xi_{1,0}=\mu^{(|\pi|)}(z^\kappa(y)),
    \qquad
    E_1=Q^{|\pi|,\kappa},
\]
and, for $j=1,\dots,|\pi|$,
\[
    \xi_{j+1}=p_{\beta_{I_j}}^\kappa(x),
    \qquad
    \xi_{j+1,0}=p_{\beta_{I_j}}^\kappa(y),
    \qquad
    E_{j+1}=P_{\beta_{I_j}}^\kappa.
\]
The hypotheses $|\xi_i-\xi_{i,0}|\le E_i$ of \cref{lem:xi-lemma} follow from
\eqref{eq:bounds-p-and-mu-b} for $i=1$, and from
\eqref{eq:bounds-p-and-mu-a} for $i=2,\dots,|\pi|+1$.  Moreover,
$|\xi_{1,0}|+E_1=M_{|\pi|}^\kappa$ and
$|\xi_{j+1,0}|+E_{j+1}=\widehat p_{\beta_{I_j}}^\kappa$.  Then \eqref{eq:xi-lemma} reads
\[
\begin{split}
\left|
\mu^{(|\pi|)}(z^\kappa(x))
\odot
\bigodot_{I\in\pi}p_{\beta_I}^\kappa(x)
-
\mu^{(|\pi|)}(z^\kappa(y))
\odot
\bigodot_{I\in\pi}p_{\beta_I}^\kappa(y)
\right|
\le
M_{|\pi|}^\kappa
\odot
\bigodot_{I\in\pi}\widehat p_{\beta_I}^\kappa
-
\left|
\mu^{(|\pi|)}(z^\kappa(y))
\odot
\bigodot_{I\in\pi}p_{\beta_I}^\kappa(y)
\right|.
\end{split}
\]
Summing over all $\pi\in\Pi_\alpha$ gives
\begin{equation}\label{eq:bound-iteration-activation-variation}
    |h_\alpha^\kappa(x)-h_\alpha^\kappa(y)|\le B_\alpha^\kappa.
\end{equation}
Differentiating $z^{\kappa+1}=W^{\kappa+1}a^\kappa+b^{\kappa+1}$ with respect to $\alpha$,
\[
    p_\alpha^{\kappa+1}(x)-p_\alpha^{\kappa+1}(y)
    =
    W^{\kappa+1}\bigl(h_\alpha^\kappa(x)-h_\alpha^\kappa(y)\bigr).
\]
Taking absolute values and using \eqref{eq:bound-iteration-activation-variation} gives \eqref{eq:affine-successor-variation-bound}. Note in particular that if $\kappa=L$, this yields
\eqref{eq:output-variation-bound} since $z^{L+1}=f$.
\end{proof}

\subsubsection{Activation variation bounds}\label{subsec:activation-variation}
The first step is to construct bounds for the activation function $\mu$ analogous to \eqref{eq:derivative-difference-bound}, these are used to construct the bounds \eqref{eq:bounds-p-and-mu-b}.

\begin{lemma}[Activation variation bounds]\label{lem:tanh-activation-bound}
Let $\mu=\tanh$ and fix a hidden layer $k$.
Fix an integer $m\ge1$ and let $N_m\in \N_0$. Define $Q^{m,k} \colon \R^w\to\R^w$ by
\begin{equation}\label{eq:tanh-activation-bound}
\begin{split}
Q^{m,k}(\rho)
&\coloneqq\sum_{r=1}^{N_m}\frac{\rho^r}{r!}\odot|\mu^{(m+r)}(z^k(y))|\\
&\quad +\rho^{N_m+1}\odot
\frac{(m+N_m+1)!}{(N_m+1)!}2^{m+N_m+2}
\exp\{-2\max(0,|z^k(y)|-\rho)\}.
\end{split}
\end{equation}
Here powers, maxima, and exponentials are applied componentwise.  If
$\rho^k\in\R_+^w$ satisfies
\begin{equation}\label{eq:rho-bound}
    |z^k(x)-z^k(y)|\leq \rho^k,
    \qquad x\in\mathcal B_{\boldsymbol{\eps}}(y),
\end{equation}
then
\begin{equation}\label{eq:activation-bounds}
    |\mu^{(m)}(z^k(x))-\mu^{(m)}(z^k(y))|\leq Q^{m,k}(\rho^k),
    \qquad x\in\mathcal B_{\boldsymbol{\eps}}(y).
\end{equation}

\end{lemma}
\begin{proof}
See \cref{app:tanh-activation-bound-proof}.
\end{proof}

\begin{remark}
The algorithms below are written for a general derivative order.  In the
numerical implementation used here, exact formulas for $\tanh^{(r)}$ are
included up to order six; for derivatives through order four we therefore
set $N_m\coloneqq6-m$, $m=1,\dots,4$.  Higher derivative orders can be handled
by supplying the corresponding higher activation derivatives and choosing $N_m$ accordingly.
\end{remark}

\subsection{Neural network derivative bound algorithms}

We now have the machinery needed to construct the bounds \eqref{eq:derivative-difference-bound}. We provide the details and pseudocode for the first order case, then recursively extend this algorithm to deal with a general multi-index $\alpha$ of arbitrary order.

\subsubsection{First order derivative bounds}
\label{subsec:first-derivative-propagation}

We give a \emph{vectorised} algorithm that handles all coordinate directions in one pass through the layers of the network and computes the activation bounds for all derivatives up to a given order.

\begin{algorithm}
\caption{Activation and first derivative bounds}\label{alg:first-deriv-bounds}
\begin{algorithmic}[1]
\Require Network $f\in\cF_{L,w}$, centre $y\in\R^n$, radius vector $\boldsymbol{\eps}\in[0,\infty)^n$, and maximum activation order $q\in\N$.
\Ensure A vector $P^{L+1}=\alg{1}(f,y,\boldsymbol{\eps})\in\R^n$ whose $\ell$-th component satisfies \eqref{eq:derivative-difference-bound} with $\alpha=e_\ell$, and activation bounds $Q^{m,k}$ for $m=1,\dots,q$ satisfying \eqref{eq:activation-bounds}.
\State Set $p_0^1\gets W^1$ and $P^1\gets0$.
\For{$k=1,\dots,L$}
    \State Set $\widehat p^k\gets |p_0^k|+P^k$ and $\rho^k\gets\widehat p^k\boldsymbol{\eps}$.
    \State Set $Q^{m,k}\gets Q^{m,k}(\rho^k)$, $m=1,\dots,q$, using \eqref{eq:tanh-activation-bound}.
    \State Set $M_1^k \gets |\mu'(z^k(y))|+Q^{1,k}$ as in \eqref{eqdef:phat-Mk}.
    \State Set $B^k\gets M_1^k\odot \widehat p^k-|\mu'(z^k(y))\odot p_0^k|$.
    \State Set $P^{k+1}\gets |W^{k+1}|B^k$.
    \State Set $p_0^{k+1}\gets W^{k+1}\bigl(\mu'(z^k(y))\odot p_0^k\bigr)$.
\EndFor
\State \Return $\alg{1}(f,y,\boldsymbol{\eps})\gets P^{L+1}$ and $\p*{Q^{m,k}}_{m, k}$.
\end{algorithmic}
\end{algorithm}

\begin{theorem}[Correctness of \cref{alg:first-deriv-bounds}]\label{thm:first-order-iteration}
Given $f\in\cF_{L,w}$, $y\in\R^n$,
$\boldsymbol{\eps}\in[0,\infty)^n$, and $q\in \N$, for every
$\ell\in\{1,\dots,n\}$, the $\ell$-th component of the returned vector
$P^{L+1}=\alg{1}(f,y,\boldsymbol{\eps})$ in
\cref{alg:first-deriv-bounds} satisfies
\eqref{eq:derivative-difference-bound} with $\alpha=e_\ell$.  The
computed $Q^{m,k}$ satisfy \eqref{eq:activation-bounds} for every
$k\in\{1,\dots,L\}$ and $m=1,\dots,q$.
\end{theorem}

\begin{proof}
First note that the chain rule for $\alpha=e_\ell$ is
\begin{equation*}
    h_\ell^k
    =\mu'(z^k)\odot p_\ell^k.
\end{equation*}
Hence, by induction on $k$, each matrix $p_0^k$ computed in \cref{alg:first-deriv-bounds} has $\ell$-th column given by $p_\ell^k(y)$.
By construction, for $\kappa=k$ and $\alpha = e_\ell$, the $\ell$-th column $B_{\ell}^k$ of the computed matrix $B^k$ is exactly $B_{\alpha}^\kappa$ defined in \eqref{eq:abstract-B-alpha}, and so the $\ell$-th column $P_\ell^{k+1}$ of the computed matrix $P^{k+1}$ is exactly the right-hand side of \eqref{eq:affine-successor-variation-bound}. For the inequality \eqref{eq:affine-successor-variation-bound} to hold at each layer, we need to verify the hypotheses \eqref{eq:bounds-p-and-mu}. We also do this by induction.

We first claim that the following implications hold for $k=1,\dots,L$,
\[
\begin{aligned}
    P_\ell^k \text{ satisfies }  \eqref{eq:bounds-p-and-mu-a} \text{ with } \kappa = k \text{ for each } \ell &\implies
    Q^{1,k} \text{ satisfies } \eqref{eq:bounds-p-and-mu-b} \text{ with } \kappa = k \\
    \qquad&\implies
    P_\ell^{k+1}
    \text{ satisfies } \eqref{eq:bounds-p-and-mu-a}
    \text{ with } \kappa=k+1 \text{ for each } \ell.
\end{aligned}
\]

For the first implication, note that the columns of $\widehat p^k$ computed in
the algorithm bound the columns of the Jacobian of $z^k$ on
$\mathcal B_{\boldsymbol{\eps}}(y)$, and
so the fundamental theorem of calculus applied to the components of $z^k$ gives
\[
\abs{z^k(x)-z^k(y)} 
\le \int_0^1 |\nabla z^k(y+\tau(x-y))|\,|x-y|\,\dd\tau
\le \widehat p^k\boldsymbol{\eps}
= \rho^k.
\]
Then, by \cref{lem:tanh-activation-bound}, the computed
$Q^{m,k}=Q^{m,k}(\rho^k)$ satisfy \eqref{eq:activation-bounds} for
$m=1,\dots,q$, and in particular they satisfy
\eqref{eq:bounds-p-and-mu-b} with $\kappa=k$ and $m=1$.

For the second implication, note that the hypotheses of
\cref{lem:bound-iteration} are satisfied by $P_\ell^k$ and $Q^{1, k}$ for
$\kappa = k$ and $\alpha=e_\ell$. Thus, since the $\ell$-th column of $P^{k+1}$ is
given by $P_\ell^{k+1} = |W^{k+1}|B_{e_\ell}^k$, it satisfies the bound
\eqref{eq:affine-successor-variation-bound} and since
$\Gamma_{e_\ell} = \{e_\ell\}$, it satisfies \eqref{eq:bounds-p-and-mu-a} for
$\kappa = k+1$.

It only remains to check the base case, i.e., that $P_\ell^1 = 0$ satisfies \eqref{eq:bounds-p-and-mu-a} for $\kappa = 1$. But this is immediate since $p_\ell^1 \equiv W^1 e_\ell$ is constant. Iterating the above implications shows that
$P_\ell^{L+1}$ satisfies \eqref{eq:bounds-p-and-mu-a} with
$\kappa=L+1$ and $\alpha=e_\ell$, which is \eqref{eq:derivative-difference-bound}
since $z^{L+1}=f$.

\end{proof}

Note that \cref{alg:first-deriv-bounds} yields a bound on the function $f$ itself, since
\begin{equation}\label{eq:zeroth-order-derivative-bound}
    \mathcal E_f^0(y,\boldsymbol{\eps})
    \coloneqq
    |f(y)|+\sum_{\ell=1}^n\eps_\ell
    \mathcal E_f^{e_\ell}(y,\boldsymbol{\eps}) 
    \geq \abs*{f(x)}, \quad x\in\mathcal B_{\boldsymbol{\eps}}(y).
\end{equation}

\subsubsection{Higher order derivative propagation}\label{subsec:higher-order-derivative-propagation}

This section extends the first order algorithm to derivatives of any order
$s\ge2$ for which the required activation bounds are available. Using \cref{lem:bound-iteration}, this is a straightforward inductive procedure.

\begin{algorithm}
\caption{Higher order derivative bound}\label{alg:higher-deriv-bounds}
\begin{algorithmic}[1]
\Require Network $f\in\cF_{L,w}$, centre $y\in\R^n$, radius vector $\boldsymbol{\eps}\in[0,\infty)^n$, and multi-index $\alpha$ with $s\coloneqq|\alpha|\ge2$.
\Ensure Scalar $\alg{s}(\alpha,f,y,\boldsymbol{\eps})$ satisfying \eqref{eq:derivative-difference-bound}.
\State Compute the activation bounds $Q^{m,k}$, $m=1,\dots,s$, $k=1,\dots,L$, using \cref{alg:first-deriv-bounds} with $q=s$.
\State Recursively compute the lower order bounds $P_\gamma^k$ satisfying
\eqref{eq:bounds-p-and-mu-a} for $\gamma\in\Gamma_\alpha\setminus\{\alpha\}$ and $k=1,\dots,L$: use \cref{alg:first-deriv-bounds} with $q=s$ when $|\gamma|=1$, and for $r=2,\dots,s-1$ compute using $\alg{r}(\gamma,f,y,\boldsymbol{\eps})$ for each $\gamma\in\Gamma_\alpha$ with $|\gamma|=r$.
\State Set $p_{\alpha,0}^1\gets0$ and $P_\alpha^1\gets0$.
\For{$k=1,\dots,L$}
    \State Set $\widehat p_\gamma^k$ for $\gamma\in\Gamma_\alpha$ and $M_m^k$, $m=1,\dots,s$, by \eqref{eqdef:phat-Mk}.
    \State Set $B_\alpha^k$ by \eqref{eq:abstract-B-alpha} with $\kappa=k$.
    \State Set $p_{\alpha,0}^{k+1}\gets W^{k+1}h_\alpha^k(y)$, using \eqref{eq:abstract-product-expansion}.
    \State Set $P_\alpha^{k+1}\gets |W^{k+1}|B_\alpha^k$.
\EndFor
\State \Return $\alg{s}(\alpha,f,y,\boldsymbol{\eps})\gets P_\alpha^{L+1}$.
\end{algorithmic}
\end{algorithm}

\begin{theorem}[Correctness of \cref{alg:higher-deriv-bounds}]\label{thm:higher-order-iteration}
Let $f\in\cF_{L,w}$, $y\in\R^n$,
$\boldsymbol{\eps}\in[0,\infty)^n$, and multi-index
$\alpha$ with $|\alpha|\ge2$ be given.
Then \cref{alg:higher-deriv-bounds} returns a number
$\alg{\abs{\alpha}}(\alpha,f,y,\boldsymbol{\eps})$ satisfying
\begin{equation}\label{eq:derivative-diff-bound-higher-ord-proof}
    |\partial^\alpha f(x)-\partial^\alpha f(y)|
    \leq \alg{\abs{\alpha}}(\alpha,f,y,\boldsymbol{\eps}),
    \qquad x\in\mathcal B_{\boldsymbol{\eps}}(y).
\end{equation}
\end{theorem}
\begin{proof}
We will argue by induction on $s=|\alpha|$. First note that the activation bounds
\eqref{eq:bounds-p-and-mu-b} for $m=1,\dots,s$ and $k=1, \dots, L$ are supplied in the first step of the algorithm, justified by \cref{thm:first-order-iteration}. It remains to show that the bound \eqref
{eq:bounds-p-and-mu-a} holds for each $\kappa \leq L$ since then \eqref{eq:derivative-diff-bound-higher-ord-proof} is equivalent to \eqref{eq:affine-successor-variation-bound}.

For $s=2$, every element of
$\Gamma_\alpha\setminus\{\alpha\}$ has order one, so
\cref{thm:first-order-iteration} gives the required lower order bounds
\eqref{eq:bounds-p-and-mu-a}.

Suppose $s\ge3$ and that for every
multi-index $\beta$ with $2\le|\beta|<s$, the recursively computed quantities
$P_\beta^k$ satisfy \eqref{eq:bounds-p-and-mu-a} for each layer $k=1,\dots,L$.
Thus, every $\gamma\in\Gamma_\alpha\setminus\{\alpha\}$ has $|\gamma|<s$, and so
the bound \eqref{eq:bounds-p-and-mu-a} is supplied either by
\cref{thm:first-order-iteration} when $|\gamma|=1$ or by the induction
hypothesis when $|\gamma|\ge2$.

With these lower order bounds available at every layer, we observe that by the bound propagation inequality \eqref{eq:affine-successor-variation-bound} of \cref{lem:bound-iteration}, and the algorithm's construction of $B_\alpha^k$ and $P_\alpha^{k+1}$, if \eqref{eq:bounds-p-and-mu-a} is satisfied for $\kappa=k < L $ by $P_\alpha^k$, then it is automatically satisfied for $\kappa = k+1$ by $P_\alpha^{k+1}$. But at the first layer, $p_\alpha^1\equiv0$, and so $P_\alpha^1=0$ satisfies \eqref{eq:bounds-p-and-mu-a} for $\kappa=1$ and hence \eqref{eq:bounds-p-and-mu} holds up to $\kappa=L$. Thus, for every $\alpha$ with $\abs{\alpha}=s$, the computed
$P_\alpha^k$ satisfy the layerwise preactivation variation bounds
\eqref{eq:bounds-p-and-mu-a} for $k=1,\dots,L$, completing the induction.

Finally, by \cref{lem:bound-iteration}, this means \eqref{eq:affine-successor-variation-bound} is satisfied for $\kappa=L$ by
$P_\alpha^{L+1}=\abs{W^{L+1}}B_\alpha^L$, or equivalently
\eqref{eq:derivative-diff-bound-higher-ord-proof} holds, since $z^{L+1}=f$.
\end{proof}

We finish the section with the simple but important observation that the computed bounds cannot blow up on a bounded domain. 

\begin{lemma}[Local boundedness of the NN bounds]
\label{lem:network-bound-local-boundedness}
For every multi-index $\alpha$, the bound
$\mathcal E_f^\alpha(y,\boldsymbol{\eps})$ defined by
\cref{eq:derivative-bound,eq:zeroth-order-derivative-bound} is
continuous in $(y,\boldsymbol{\eps})\in\R^n\times[0,\infty)^n$.
Consequently, for every compact $K\subset\R^n$ and every
$\overline{\boldsymbol{\eps}}\in[0,\infty)^n$,
\[
 \sup_{y\in K,\ 0\leq\boldsymbol{\eps}\leq
 \overline{\boldsymbol{\eps}}}
 \mathcal E_f^\alpha(y,\boldsymbol{\eps})<\infty.
\]
\end{lemma}

\begin{proof}
The NN and each of its derivatives are continuous. 
\cref{alg:first-deriv-bounds,alg:higher-deriv-bounds} use finitely many
evaluations of these functions and finitely many additions, products,
absolute values, maxima, exponentials, and matrix multiplications. Thus their outputs are continuous. The claim for $\abs{\alpha} > 0$ follows from the extreme value theorem, and the assertion
for $\alpha=0$ follows from
\eqref{eq:zeroth-order-derivative-bound}. 
\end{proof}

\section{Residual-based verification for heat and wave equations}\label{sec:applications}

In this section, we return to a posteriori verification of neural network-based PDE solvers. We first provide concrete residual energy estimates of the form of \cref{eq:heat-h1-framework,eq:wave-h1-framework}, and discuss how these residual norms will be rigorously bounded using the quadrature framework described in \cref{subsec:quadrature}. We then detail how the derivative bounds from \cref{sec:activation} are used to compute the local error bound \eqref{eq:cell-model-errors} for each residual norm.

\subsection{Explicit energy estimates}\label{sec:energy-estimates}

We first derive explicit energy estimates for the heat and wave equations, which will be used to establish the residual framework \eqref{eq:heat-wave-model-framework}.

\begin{notation}\label{not:heat-wave-operator}
Throughout, \(H_0^1(\Omega)\) is equipped with the norm
\[
    \norm{v}_{H_0^1(\Omega)}\coloneqq\norm{\nabla v}_{L^2(\Omega)},
\]
and \(H^{-1}(\Omega)\) denotes the dual space with respect to this norm. In the Sobolev space $H^2(\Omega)$, we use the norm
\[
    \norm{w}_{H^2(\Omega)}
    :=
    \left(
        \norm{w}_{L^2(\Omega)}^2
        +\norm{\nabla w}_{L^2(\Omega)}^2
        +\norm{D^2w}_{L^2(\Omega)}^2
    \right)^{1/2}.
\]
Let \(A:H_0^1(\Omega)\to H^{-1}(\Omega)\) be the operator defined by
\[
    \langle Av,w\rangle_{H^{-1}(\Omega),H_0^1(\Omega)}
    :=a(v,w)
    :=\int_\Omega \nabla v\cdot\nabla w\,dx,
    \qquad v,w\in H_0^1(\Omega).
\]
We use the same symbol \(A\) for its \(L^2(\Omega)\)-realisation, whose domain is
\[
    D(A):=\{v\in H_0^1(\Omega):Av\in L^2(\Omega)\}.
\]
The restriction \(A:D(A)\subset L^2(\Omega)\to L^2(\Omega)\) is the positive
self-adjoint Dirichlet Laplacian \(A=-\Delta\).
Let \((w_k,\lambda_k)_{k\ge1}\) denote
the eigenpairs of \(A\) that form an orthonormal basis of
\(L^2(\Omega)\) and are orthogonal in \(H_0^1(\Omega)\)
\citep[Section 6.5, Theorem 1]{evans_2010}, with
\begin{equation}\label{eq:dirichlet-eigenbasis}
    Aw_k=\lambda_k w_k,
    \qquad
    0<\lambda_1\le\lambda_2\le\cdots .
\end{equation}
For \(m\ge1\), write
\(
    V_m\coloneqq \operatorname{span}\{w_1,\dots,w_m\},
\)
and, whenever Galerkin approximations are used, take \(u_m(t)\in V_m\) to
satisfy the corresponding weak form for all test functions in \(V_m\) and a.e.
\(t\in(0,T)\).
\end{notation}

We then have the explicit Poincar\'e inequality
\begin{equation}\label{eq:poincare-h01}
    \norm{v}_{L^2(\Omega)}
    \le
    \lambda_1^{-1/2}\norm{v}_{H_0^1(\Omega)},
    \qquad v\in H_0^1(\Omega).
\end{equation}

\begin{lemma}\label{lem:convex-domain-A-H2}
Assume that \(\Omega\) is convex. Then every \(w\in D(A)\) belongs
to \(H^2(\Omega)\cap H_0^1(\Omega)\), and
\[
    \norm{w}_{H^2(\Omega)}
    \le
    C_\Omega\norm{Aw}_{L^2(\Omega)},
    \qquad
    C_\Omega
    =
    \left(
        1+\frac1{\lambda_1}+\frac1{\lambda_1^2}
    \right)^{1/2}.
\]
\end{lemma}

\begin{proof}
Since \(w\in D(A)\), we have \(Aw\in L^2(\Omega)\) and \(w\) solves the
Dirichlet problem \(-\Delta w=Aw\) in the distributional sense.  Convex-domain elliptic regularity
\citep[Theorem 3.2.1.2]{grisvard_1985} gives
\(w\in H^2(\Omega)\cap H_0^1(\Omega)\).  Then the classical Miranda-Talenti estimate gives  
\[
    \norm{D^2 w}_{L^2(\Omega)}
    \le
    \norm{Aw}_{L^2(\Omega)}.
\]
Moreover, integration by parts and Poincar\'e's inequality give
\[
    \norm{\nabla w}_{L^2}
    \le
    \lambda_1^{-1/2}\norm{Aw}_{L^2},
    \qquad
    \norm{w}_{L^2}
    \le
    \lambda_1^{-1}\norm{Aw}_{L^2}.
\]
Adding these three inequalities gives the desired result.
\end{proof}

\subsubsection{General heat equation energy estimates}

\begin{theorem}\label{thm:heat-estimates}
Let \(\Omega\) be an open, bounded set and $\kappa, T > 0$.
Consider the heat equation
\begin{equation}\label{eq:heat-general}
\begin{cases}
\partial_t u-\kappa\Delta u=F & \text{in }\Omega\times(0,T],\\
u=0 & \text{on }\partial\Omega\times[0,T],\\
u=g & \text{on }\Omega\times\set{t=0},
\end{cases}
\end{equation}
where \(F\in L^2(0,T;L^2(\Omega))\) and \(g\in L^2(\Omega)\).  There is a
unique weak solution
\[
    u\in C([0,T];L^2(\Omega))
    \qquad
    u'\in L^2(0,T;H^{-1}(\Omega)),
\]
for which the following is true:
\begin{enumerate}[label=\textup{(\arabic*)}]
\item\label{thm:heat-estimates-l2} $u \in L^2(0,T;H_0^1(\Omega))$ with the bound
\begin{equation}\label{eq:explicit-L2-estimate}
\begin{split}
    &\sup_{0\le t\le T}\norm{u(t)}_{L^2(\Omega)}
    + \norm{u}_{L^2(0,T;H_0^1(\Omega))}
    + \norm{u'}_{L^2(0,T;H^{-1}(\Omega))}
    \\
    &\qquad\le
    \alpha^{\mathrm H}\norm{g}_{L^2(\Omega)}
    +
    \beta^{\mathrm H}
    \norm{F}_{L^2(0,T;L^2(\Omega))},
\end{split}
\end{equation}
where
\begin{equation}\label{eq:heat-alpha-beta-constants}
    \alpha^{\mathrm H}\coloneqq 1+\frac1{\sqrt{\kappa}}+\sqrt{\kappa},
    \qquad
    \beta^{\mathrm H}\coloneqq \frac1{\sqrt{\lambda_1}}
    \left(
        2+
        \frac1{\sqrt{\kappa}}
        +
        \frac1{\kappa}
    \right).
\end{equation}

\item\label{thm:heat-estimates-h1} If \(\Omega\) is convex and
\(
    g\in H_0^1(\Omega),
\)
then
\begin{equation}\label{eq:u-increased-reg}
    u\in L^\infty(0,T;H_0^1(\Omega))
    \cap L^2(0,T;H^2(\Omega)),
    \qquad
    u'\in L^2(0,T;L^2(\Omega)).
\end{equation}
Moreover,
\begin{equation}\label{eq:explicit-H1-H2-estimate}
\begin{split}
    &\operatorname*{ess\,sup}_{0\le t\le T}\norm{u(t)}_{H_0^1(\Omega)}
    +
    \norm{u}_{L^2(0,T;H^2(\Omega))}
    +
    \norm{u'}_{L^2(0,T;L^2(\Omega))}\\
    &\qquad\le
    C_1
    \norm{g}_{H_0^1(\Omega)}
    +
    C_2
    \norm{F}_{L^2(0,T;L^2(\Omega))},
\end{split}
\end{equation}
where \(C_\Omega\) is the constant from \cref{lem:convex-domain-A-H2}, and
\begin{equation}\label{eq:C-H1-H2}
\begin{gathered}
    C_1
    =
    1+\sqrt{\kappa}+\frac{C_\Omega}{\sqrt{\kappa}},
    \quad
    C_2
    =
    2+\frac1{\sqrt{\kappa}}+\frac{C_\Omega}{\kappa}.
\end{gathered}
\end{equation}
\end{enumerate}
\end{theorem}

\begin{proof}[Proof of \cref{thm:heat-estimates}]
The existence and uniqueness statements are standard results \citep[Section 7.1.2]{evans_2010}. For the estimates, we base the proofs on \citet{evans_2010}, but we keep track of the constants to obtain explicit bounds and avoid the smooth boundary assumption using \cref{lem:convex-domain-A-H2}.

\textbf{Part (1)}.
Test the $m$-th Galerkin equation for \eqref{eq:heat-general} against the corresponding approximation \(u_m(t)\). Poincar\'e's and Young's
inequalities give
\begin{equation}\label{eq:proof-L2-H1-differential}
    \frac12\frac{d}{dt}\norm{u_m(t)}_{L^2}^2
    +
    \kappa\norm{\nabla u_m(t)}_{L^2}^2
    =
    (F(t),u_m(t))_{L^2}
    \le
    \frac{\kappa}{2}\norm{\nabla u_m(t)}_{L^2}^2
    +
    \frac{1}{2\kappa\lambda_1}\norm{F(t)}_{L^2}^2.
\end{equation}
Integrating over \((0,t)\), taking the supremum over \((0,T)\), and
passing to the weak limit gives, by lower semicontinuity,
\begin{equation}\label{eq:proof-L2-energy-bounds}
\begin{split}
    \sup_{0\le t\le T}\norm{u(t)}_{L^2}
    &\le
    \norm{g}_{L^2}
    +
    \frac1{\sqrt{\kappa\lambda_1}}
    \norm{F}_{L^2(0,T;L^2)},
    \\
    \norm{u}_{L^2(0,T;H_0^1)}
    &\le
    \frac1{\sqrt{\kappa}}\norm{g}_{L^2}
    +
    \frac1{\kappa\sqrt{\lambda_1}}
    \norm{F}_{L^2(0,T;L^2)}.
\end{split}
\end{equation}
To bound the third norm, note that
the dual norm satisfies
\begin{equation}\label{eq:proof-Hminus1-dual-bounds}
    \norm{F(t)}_{H^{-1}}
    \le
    \lambda_1^{-1/2}\norm{F(t)}_{L^2},
    \qquad
    \norm{Au(t)}_{H^{-1}}
    =
    \norm{u(t)}_{H_0^1}.
\end{equation}
Thus, since $u'=F-\kappa Au$ in $L^2(0,T;H^{-1}(\Omega))$, the second inequality in \eqref{eq:proof-L2-energy-bounds} gives
\begin{equation}\label{eq:proof-uprime-Hminus1}
\begin{split}
    \norm{u'}_{L^2(0,T;H^{-1})}
    \le
    \lambda_1^{-1/2}
    \norm{F}_{L^2(0,T;L^2)}
    +
    \kappa
    \norm{u}_{L^2(0,T;H_0^1)}
    \le
    \sqrt{\kappa}\norm{g}_{L^2}
    +
    \frac2{\sqrt{\lambda_1}}
    \norm{F}_{L^2(0,T;L^2)}.
\end{split}
\end{equation}
Combining \eqref{eq:proof-L2-energy-bounds} and \eqref{eq:proof-uprime-Hminus1} proves
\eqref{eq:explicit-L2-estimate}.

\textbf{Part \textup{(2)}}.
Testing the $m$-th Galerkin equation with \(Au_m(t)\) and applying Young's
inequality gives
\begin{equation}\label{eq:proof-H1-Au-identity}
    \frac12\frac{d}{dt}\norm{\nabla u_m(t)}_{L^2}^2
    +
    \kappa\norm{Au_m(t)}_{L^2}^2
    =
    (F(t),Au_m(t))_{L^2}
    \le
    \frac1{2\kappa}\norm{F(t)}_{L^2}^2
    +
    \frac{\kappa}{2}\norm{Au_m(t)}_{L^2}^2.
\end{equation}
Integrating over \((0,t)\), taking the essential supremum over \((0,T)\), and
passing to the weak limit,
\begin{equation}\label{eq:proof-Linf-H01}
    \operatorname*{ess\,sup}_{0\le t\le T}\norm{u(t)}_{H_0^1(\Omega)}
    \le
    \norm{g}_{H_0^1(\Omega)}
    +
    \frac1{\sqrt{\kappa}}
    \norm{F}_{L^2(0,T;L^2(\Omega))},
\end{equation}
and
\begin{equation}\label{eq:proof-Au-L2}
    \norm{Au}_{L^2(0,T;L^2(\Omega))}
    \le
    \frac1{\sqrt{\kappa}}\norm{g}_{H_0^1(\Omega)}
    +
    \frac1{\kappa}
    \norm{F}_{L^2(0,T;L^2(\Omega))}.
\end{equation}
Thus, since \(\Omega\) is convex, \cref{lem:convex-domain-A-H2} implies
\begin{equation}\label{eq:proof-H2-via-Au}
    \norm{u}_{L^2(0,T;H^2(\Omega))}
    \le
    C_\Omega\norm{Au}_{L^2(0,T;L^2(\Omega))}
    \le
    C_\Omega
    \left(
        \frac1{\sqrt{\kappa}}\norm{g}_{H_0^1(\Omega)}
        +
        \frac1{\kappa}
        \norm{F}_{L^2(0,T;L^2(\Omega))}
    \right).
\end{equation}
Finally, since $u'=F-\kappa Au$ in $L^2(0,T;L^2(\Omega))$, 
\begin{equation}\label{eq:proof-uprime-L2}
    \norm{u'}_{L^2(0,T;L^2(\Omega))}
    \le
    \sqrt{\kappa}\norm{g}_{H_0^1(\Omega)}
    +
    2\norm{F}_{L^2(0,T;L^2(\Omega))}.
\end{equation}
Adding \eqref{eq:proof-Linf-H01}, \eqref{eq:proof-H2-via-Au}, and
\eqref{eq:proof-uprime-L2} proves \eqref{eq:explicit-H1-H2-estimate}.
\end{proof}

\subsubsection{General wave equation energy estimates}

For the wave equation we provide a single energy estimate matching
\hyperref[thm:heat-estimates-h1]{\cref*{thm:heat-estimates}, part~\textup{(2)}}.

\begin{theorem}\label{thm:wave-estimates}
Let \(\Omega\) be an open, bounded domain and $c, T > 0$. Consider the wave equation
\begin{equation}\label{eq:wave-general}
\begin{cases}
\partial_t^2 u-c^2\Delta u=F & \text{in }\Omega\times(0,T],\\
u=0 & \text{on }\partial\Omega\times[0,T],\\
u=g,\quad \partial_t u=h & \text{on }\Omega\times\set{t=0},
\end{cases}
\end{equation}
where
\[
    F\in L^2(0,T;L^2(\Omega)),
    \qquad
    g\in H_0^1(\Omega),
    \qquad
    h\in L^2(\Omega).
\]
There is a unique weak solution
\[
    u\in L^\infty(0,T;H_0^1(\Omega)),
    \qquad
    u'\in L^\infty(0,T;L^2(\Omega)),
    \qquad
    u''\in L^2(0,T;H^{-1}(\Omega)),
\]
for which the following estimate holds:
\begin{equation}\label{eq:explicit-wave-energy-estimate}
\begin{split}
    &\operatorname*{ess\,sup}_{0<t<T}
    \left(
        \norm{u(t)}_{H_0^1(\Omega)}
        +
        \norm{u'(t)}_{L^2(\Omega)}
    \right)
    +
    \norm{u''}_{L^2(0,T;H^{-1}(\Omega))}
    \\
    &\qquad\le
    \alpha^{\mathrm W}
    \norm{g}_{H_0^1(\Omega)}
    +
    \eta^{\mathrm W}
    \norm{h}_{L^2(\Omega)}
    +
    \beta^{\mathrm W}
    \norm{F}_{L^2(0,T;L^2(\Omega))},
\end{split}
\end{equation}
where
\begin{equation}\label{eq:wave-alpha-beta-constants}
    \alpha^{\mathrm W}\coloneqq 1+c+c^2\sqrt{T},
    \qquad
    \eta^{\mathrm W}\coloneqq 1+\frac1c+c\sqrt{T},
    \qquad
    \beta^{\mathrm W}\coloneqq \frac1{\sqrt{\lambda_1}} + cT
    +\sqrt{T}\left(1+\frac1c\right).
\end{equation}

\end{theorem}

\begin{proof}[Proof of \cref{thm:wave-estimates}]
The existence and uniqueness statements are standard; see
\citet[Section 7.2.2]{evans_2010}. It only remains to keep track of constants in the estimates. Testing the $m$-th Galerkin equation with
\(u_m'(t)\) gives
\[
    \frac12\frac{d}{dt}
    \left(
        E_m(t)^2
    \right)
    =
    (F(t),u_m'(t))_{L^2},
\]
where
\[
    E_m(t)^2
    \coloneqq
    \norm{u_m'(t)}_{L^2}^2
    +
    c^2\norm{u_m(t)}_{H_0^1}^2 .
\]
Since \(E_m\) is absolutely continuous, dividing by \(E_m(t)\) where
\(E_m(t)>0\) and applying Cauchy--Schwarz gives, for a.e. \(t\in(0,T)\),
\[
    E_m'(t)
    \le
    \norm{F(t)}_{L^2}.
\]
Moreover, \(E_m'(t)=0\) for a.e.\ \(t\) such that \(E_m(t)=0\).
Thus,
\[
    E_m(t)\le E_m(0)+\int_0^t\norm{F(s)}_{L^2}\,ds .
\]

Taking the essential supremum over \((0,T)\), applying Cauchy--Schwarz to the
resulting integral involving $F$, and passing to the weak limit gives
\begin{equation}\label{eq:proof-wave-energy-bounds}
\begin{split}
    \operatorname*{ess\,sup}_{0<t<T}\norm{u(t)}_{H_0^1}
    &\le
    \norm{g}_{H_0^1}
    +
    \frac1c\norm{h}_{L^2}
    +
    \frac{\sqrt{T}}{c}
    \norm{F}_{L^2(0,T;L^2)},
    \\
    \operatorname*{ess\,sup}_{0<t<T}\norm{u'(t)}_{L^2}
    &\le
    c\norm{g}_{H_0^1}
    +
    \norm{h}_{L^2}
    +
    \sqrt{T}
    \norm{F}_{L^2(0,T;L^2)}.
\end{split}
\end{equation}
To bound the third norm, recall \eqref{eq:proof-Hminus1-dual-bounds}.
Thus, since $u''=F-c^2Au$ in $L^2(0,T;H^{-1}(\Omega))$, the first inequality in
\eqref{eq:proof-wave-energy-bounds} gives
\begin{equation}\label{eq:proof-wave-utt}
\begin{split}
    \norm{u''}_{L^2(0,T;H^{-1})}
    &\le
    \lambda_1^{-1/2}
    \norm{F}_{L^2(0,T;L^2)}
    +
    c^2
    \norm{u}_{L^2(0,T;H_0^1)}
    \\
    &\le
    c^2\sqrt{T}\norm{g}_{H_0^1}
    +
    c\sqrt{T}\norm{h}_{L^2}
    +
    \left(
        \frac1{\sqrt{\lambda_1}}
        +
        cT
    \right)
    \norm{F}_{L^2(0,T;L^2)}.
\end{split}
\end{equation}
Adding the inequalities in \eqref{eq:proof-wave-energy-bounds} and
\eqref{eq:proof-wave-utt} proves \eqref{eq:explicit-wave-energy-estimate}.
\end{proof}

\subsubsection{Residual estimates}
Let $u_{\mathrm H}$ be the weak solution of the heat problem in
\eqref{eq:heat-wave-model-framework}, let $u_{\mathrm W}$ be the weak solution
of the wave problem in \eqref{eq:heat-wave-model-framework}, and let
$v = Bf$ be as in \cref{assumption:NN-bounds}.
Define
$e_{\mathrm H}\coloneqq u_{\mathrm H}-v$, $e_{\mathrm W}\coloneqq u_{\mathrm W}-v$,
and $v_0, \dot v_0$ as in \eqref{eq:initial-data}. Denote the PDE residuals as
\[
    R_{\mathrm H}(v)\coloneqq \partial_t v-\kappa\Delta v,
    \qquad
    R_{\mathrm W}(v)\coloneqq \partial_t^2v-c^2\Delta v.
\]
Then \(e_{\mathrm H}\) and \(e_{\mathrm W}\) are weak solutions of the heat and
wave \emph{residual} equations, respectively,
\begin{equation}\label{eq:heat-wave-residual}
\begin{array}{@{}c@{\qquad}c@{}}
\begin{cases}
\partial_t e-\kappa\Delta e=-R_{\mathrm H}(v) & \text{in }\Omega\times(0,T],\\
e=0 & \text{on }\partial\Omega\times[0,T],\\
e=g-v_0 & \text{on }\Omega\times\set{t=0},
\end{cases}
&
\begin{cases}
\partial_t^2 e-c^2\Delta e=-R_{\mathrm W}(v) & \text{in }\Omega\times(0,T],\\
e=0 & \text{on }\partial\Omega\times[0,T],\\
e=g-v_0,\quad \partial_t e=h-\dot v_0 & \text{on }\Omega\times\set{t=0}.
\end{cases}
\end{array}
\end{equation}

Then the energy estimates from \cref{sec:energy-estimates} give the following corollaries, which are the precise forms of the residual frameworks \eqref{eq:heat-h1-framework} and \eqref{eq:wave-h1-framework} for the heat and wave equations.

\begin{corollary}[Heat residual estimates]\label{cor:heat-residual-estimates}
With $e_{\mathrm H}$ defined as above, under the assumptions of
\hyperref[thm:heat-estimates-l2]{\cref*{thm:heat-estimates}, part~\textup{(1)}} and
\hyperref[thm:heat-estimates-h1]{\cref*{thm:heat-estimates}, part~\textup{(2)}},
respectively, we have
\begin{subequations}\label{eq:heat-residual-estimates}
\begin{align}
    \sup_{0\le t\le T}\norm{e_{\mathrm H}(t)}_{L^2}
    +\norm{e_{\mathrm H}}_{L^2(0,T;H_0^1)}
    +\norm{\partial_t e_{\mathrm H}}_{L^2(0,T;H^{-1})}
    &\le
    \alpha^{\mathrm H}\norm{g-v_0}_{L^2}
    +
    \beta^{\mathrm H}\norm{R_{\mathrm H}(v)}_{L^2(0,T;L^2)}, \label{eq:heat-residual-estimate-l2}
    \\
    \operatorname*{ess\,sup}_{0<t<T}\norm{e_{\mathrm H}(t)}_{H_0^1}
    +\norm{e_{\mathrm H}}_{L^2(0,T;H^2)}
    +\norm{\partial_t e_{\mathrm H}}_{L^2(0,T;L^2)}
    &\le
    C_1\norm{g-v_0}_{H_0^1}
    +
    C_2\norm{R_{\mathrm H}(v)}_{L^2(0,T;L^2)} . \label{eq:heat-residual-estimate-h1}
\end{align}
\end{subequations}
We denote the left- and right-hand sides of
\eqref{eq:heat-residual-estimate-l2} by
$\mathcal N_{\mathrm H_0}(e_{\mathrm H})$ and
$\mathcal R_{\mathrm H_0}(v)$, respectively, and those of
\eqref{eq:heat-residual-estimate-h1} by
$\mathcal N_{\mathrm H_1}(e_{\mathrm H})$ and
$\mathcal R_{\mathrm H_1}(v)$, respectively.
\end{corollary}

\begin{corollary}[Wave residual estimate]\label{cor:wave-residual-estimate}
With $e_{\mathrm W}$ defined as above, under the assumptions of \cref{thm:wave-estimates},
\begin{equation}\label{eq:wave-residual-estimate}
\begin{split}
    &\operatorname*{ess\,sup}_{0<t<T}
    \left(
        \norm{e_{\mathrm W}(t)}_{H_0^1}
        +
        \norm{\partial_t e_{\mathrm W}(t)}_{L^2}
    \right)
    +
    \norm{\partial_t^2 e_{\mathrm W}}_{L^2(0,T;H^{-1})}
    \\
    &\qquad\le
    \alpha^{\mathrm W}\norm{g-v_0}_{H_0^1}
    +
    \eta^{\mathrm W}\norm{h-\dot v_0}_{L^2}
    +
    \beta^{\mathrm W}\norm{R_{\mathrm W}(v)}_{L^2(0,T;L^2)} .
\end{split}
\end{equation}
We denote the left- and right-hand sides of
\eqref{eq:wave-residual-estimate} by
$\mathcal N_{\mathrm W}(e_{\mathrm W})$ and
$\mathcal R_{\mathrm W}(v)$, respectively.
\end{corollary}

Thus in the framework of \cref{eq:XYZ-residuals} the $X$ norm is given by one of the left-hand sides in \eqref{eq:heat-residual-estimates} or \eqref{eq:wave-residual-estimate}, and the $Y$ norm is given by the corresponding right-hand side. 

\subsection{Rigorous quadrature for residual norms}\label{sec:residual-quad}
The preceding subsection shows that the quadrature methods defined in \cref{subsec:quadrature} must be applied to the residual norms:
\begin{equation}\label{eq:quadrature-target-residual-norms}
\begin{gathered}
    \norm{g-v_0}_{L^2(\Omega)},
    \quad
    \norm{\nabla(g-v_0)}_{L^2(\Omega)},
    \quad
    \norm{h-\dot v_0}_{L^2(\Omega)},
    \\
    \norm{R_{\mathrm H}(v)}_{L^2(0,T;L^2(\Omega))},
    \quad
    \norm{R_{\mathrm W}(v)}_{L^2(0,T;L^2(\Omega))}.
\end{gathered}
\end{equation}

For concreteness, assume that \(\Omega=(0,1)^d\), $d\in\N$, and that we are given a spatial radius vector $\boldsymbol{\eps}_x\in(0,\infty)^d$ and a time radius $\eps_t>0$ satisfying \cref{assumption:rectangular-domain}. Set $\boldsymbol{\eps}\coloneqq(\boldsymbol{\eps}_x,\eps_t)$. In particular, we have the partitions
 \(
    \mathcal P_\Omega
    \coloneqq
    \set*{\mathcal B_{\boldsymbol{\eps}_x}(y_k)}_k,
  \) and
  \(  \mathcal P
    \coloneqq
    \set*{\mathcal B_{\boldsymbol{\eps}}(y_k,t_m)}_{k,m}
  \)
  of $\Omega$ and $\Omega_T$, respectively.
In this case, the boundary function has the explicit form
\[
B(x) = \prod_{i=1}^d x_i(1-x_i),
\]
whose explicit moduli $\omega_B^\beta$ satisfying \cref{ass:B-moduli} are recorded in
\cref{app:boundary-factor-moduli}. Since the residual estimates and their bounds are independent of training methodology, we need only assume that the approximation is given by $v = Bf$ with $f \in \cF_{L, w}$ for some fixed $L, w \in \N$. 

With the setup of the quadrature methods in \cref{subsec:quadrature}, the only work remaining is to construct computable cellwise bounds \eqref{eq:cell-derivative-coefficients} for each integrand in \eqref{eq:quadrature-target-residual-norms}. The NN derivative
bounds entering these quantities are those computed by
\cref{alg:first-deriv-bounds,alg:higher-deriv-bounds} and will be referenced with the notation $\mathcal E_f^\alpha$ defined in \cref{eq:derivative-bound}, except for the special case $\alpha=0$  which is given by \cref{eq:zeroth-order-derivative-bound}.

\subsubsection{Initial displacement residuals}\label{sec:initial-displacement-quad}

Set $e_0\coloneqq g-v_0=g-Bf(\cdot,0)$.  For the initial data residuals, let
$C_k=\mathcal B_{\boldsymbol{\eps}_x}(y_k)\in\mathcal P_\Omega$ and set
$\boldsymbol{\eps}_0\coloneqq(\boldsymbol{\eps}_x,0)$ and
$h_x\coloneqq\norm{\boldsymbol{\eps}_x}_\infty$.  Let
$e_1,\dots,e_{d+1}$ denote the spacetime coordinate multi-indices, and
identify spatial multi-indices with spacetime multi-indices having zero time
component.  
For multi-indices $\beta \in \N_0^{d}$ and $\gamma \in \N_0^{d+1}$, we formally define the quantities
\begin{equation}\label{eq:initial-factor-bounds}
     \mathsf G_k^\beta
 \coloneqq
 \mathcal D_g^\beta(y_k,\boldsymbol{\eps}_x),
 \qquad
 \mathsf B_k^\beta
 \coloneqq |\partial^\beta B(y_k)|+\omega_B^\beta(h_x),
 \qquad
 \mathsf F_k^\gamma
 \coloneqq
 \mathcal E_f^\gamma((y_k,0),\boldsymbol{\eps}_0).
\end{equation}

By \cref{thm:first-order-iteration,thm:higher-order-iteration}, $\mathsf F^\gamma_k$ is well-defined and computable for any order of $\gamma$, and, as mentioned above, $\mathsf{B}^\beta_k$ are explicitly available via \cref{app:boundary-factor-moduli} for $\abs{\beta}\leq 4$. For $\mathsf{G}^\beta_k$, \cref{ass:data-derivative-bounds} requires $\abs{\beta}\leq 1$ to apply $Q^0$ and $\abs{\beta}\leq 2$ to apply $Q^1$.

For every spatial multi-index $\alpha$ with
$|\alpha|\leq 4$, we can use Leibniz's rule to compute bounds
\begin{equation}\label{eq:initial-product-bound}
 \mathsf V_k^\alpha
 \coloneqq
 \sum_{\beta\leq\alpha}\binom{\alpha}{\beta}
 \mathsf B_k^\beta\mathsf F_k^{\alpha-\beta}
 \geq\sup_{x\in C_k}|\partial^\alpha(Bf)(x,0)|.
\end{equation}
Thus, for the quadrature rules $Q^0$ and $Q^1$, respectively, under their relevant assumptions, we compute the required bounds \eqref{eq:cell-derivative-coefficients} as
\begin{equation}\label{eq:initial-displacement-coefficients}
\begin{aligned}
 b_{q,C_k}(e_0)
 &\coloneqq \mathsf G_k^{e_q}+\mathsf V_k^{e_q}
 \geq\sup_{x\in C_k}|\partial_qe_0(x)|,
 \\
 h_{q\ell,C_k}(e_0)
 &\coloneqq
 \mathsf G_k^{e_q+e_\ell}
 +\mathsf V_k^{e_q+e_\ell}
 \geq\sup_{x\in C_k}|\partial_{q\ell}e_0(x)|.
\end{aligned}
\end{equation}

Given these bounds, each quadrature estimate and error bound is computed for $\phi=e_0$ exactly as written in \cref{subsec:quadrature}. 

\subsubsection{Initial displacement gradient residual}\label{sec:initial-displacement-gradient-quad}

As indicated by \cref{rem:vector-residual-bounds}, we consider the components
of $\nabla e_0$ by applying the construction and regularity assumptions above
to each $\partial_j e_0$. This yields bounds for $j,q,\ell \in \set{1,\dots,d}$
which apply to the rules $Q^0$ and $Q^1$ respectively:
\begin{equation}\label{eq:initial-gradient-coefficients}
\begin{aligned}
    b_{q,C_k}(\partial_j e_0)
&\coloneqq h_{qj,C_k}(e_0)
\geq
\sup_{x\in C_k}|\partial_q(\partial_j e_0)(x)|,
 \\
 h_{q\ell,C_k}(\partial_j e_0)
 &\coloneqq
 \mathsf G_k^{e_j+e_q+e_\ell}
 +\mathsf V_k^{e_j+e_q+e_\ell}
 \geq\sup_{x\in C_k}|\partial_{q\ell}(\partial_j e_0)(x)|.
\end{aligned}
\end{equation}
This yields componentwise bounds 
$\mathcal U_{\mathcal P_\Omega}^r(\partial_j e_0)$, and hence the global bound $\mathcal U_{\mathcal P_\Omega}^r(\nabla e_0)$ as defined in \cref{rem:vector-residual-bounds}. 

\subsubsection{Initial velocity residual} \label{sec:initial-velocity-quad}

For the initial velocity residual, we consider
\[
    e_{\mathrm v}\coloneqq h-\dot v_0
    =h-B\partial_t f(\cdot,0).
\]
Set $\tau=e_{d+1}$.  The corresponding data and product bounds are, for
spatial multi-indices $\alpha$ with $|\alpha|\leq2$,
\begin{equation}\label{eq:initial-velocity-factor-bounds}
\begin{aligned}
    \mathsf H_k^\alpha
    &\coloneqq
    \mathcal D_h^\alpha(y_k,\boldsymbol{\eps}_x),
    &
    \mathsf V_{\mathrm v,k}^\alpha
    &\coloneqq
    \sum_{\beta\leq\alpha}\binom{\alpha}{\beta}
    \mathsf B_k^\beta\mathsf F_k^{\alpha-\beta+\tau}
    \geq
    \sup_{x\in C_k}
    |\partial^\alpha(B\partial_t f)(x,0)|.
\end{aligned}
\end{equation}
Here $\mathsf H_k^\alpha$ bounds
$\sup_{x\in C_k}|\partial^\alpha h(x)|$, and the sum is only over spatial indices. Thus, for the quadrature rules
$Q^0$ and $Q^1$
respectively, under the relevant assumptions that allow us to compute \eqref{eq:initial-velocity-factor-bounds}, the required bounds
\eqref{eq:cell-derivative-coefficients} are given by
\begin{equation}\label{eq:initial-velocity-coefficients}
\begin{aligned}
 b_{q,C_k}(e_{\mathrm v})
 &\coloneqq
 \mathsf H_k^{e_q}+\mathsf V_{\mathrm v,k}^{e_q}
 \geq
\sup_{x\in C_k}|\partial_qe_{\mathrm v}(x)|,
 \\
 h_{q\ell,C_k}(e_{\mathrm v})
 &\coloneqq
 \mathsf H_k^{e_q+e_\ell}
 +\mathsf V_{\mathrm v,k}^{e_q+e_\ell}
 \geq
\sup_{x\in C_k}|\partial_{q\ell}e_{\mathrm v}(x)|.
\end{aligned}
\end{equation}
Given these bounds, each quadrature estimate and error bound is computed for
$\phi=e_{\mathrm v}$ exactly as written in \cref{subsec:quadrature}. 

\subsubsection{PDE residuals}\label{sec:pde-residual-quad}

The PDE residuals are of the form $R_{\mathfrak p}(v)$, where $\mathfrak p\in\{\mathrm H,\mathrm W\}$. They involve space and time and include second derivatives, so they present the most challenging quadrature problem.
Since the integrand is smooth, we use the affine rule $Q^1$ on a spacetime partition $\mathcal P$ with anisotropic cells
$C=\zeta_C+\prod_{q=1}^{d+1}[-\eps_{C,q},\eps_{C,q}]$, and allow adaptive
refinement.  We construct $h_{q\ell,C}^{\mathfrak p}$ satisfying
\eqref{eq:cell-derivative-coefficients}, i.e., 
\begin{equation}\label{eq:pde-residual-second-coefficients}
 h_{q\ell,C}^{\mathfrak p}
 \geq\sup_{\zeta\in C}|\partial_{q\ell}R_{\mathfrak p}(v)(\zeta)|,
 \qquad q,\ell=1,\dots,d+1,
\end{equation}
in \hyperref[app:initial-data-remainders]{Appendix~\ref*{app:initial-data-remainders}}, using network derivative bounds through
order four.

\subsection{Verification algorithms}
\label{subsec:verification-algorithms}

We give an informal \emph{local} verification algorithm and a \emph{global} one
that together affirmatively answer \cref{question-sequence} and illustrate the residual verification framework for our model problems.
Adding the residual norm bounds gives computable bounds for the quadrature rules $r\in\{0,1\}$, which we denote by
\begin{equation}\label{eq:computable-residual-bounds}
\begin{aligned}
    \widehat{\mathcal R}_{\mathrm H_0}^{r}
    (v;\mathcal P_\Omega,\mathcal P)
    &\coloneqq
    \alpha^{\mathrm H}\mathcal U_{\mathcal P_\Omega}^r(e_0)
    +\beta^{\mathrm H}\mathcal U_{\mathcal P}^1(R_{\mathrm H}(v)),
    \\
    \widehat{\mathcal R}_{\mathrm H_1}^{r}
    (v;\mathcal P_\Omega,\mathcal P)
    &\coloneqq
    C_1\mathcal U_{\mathcal P_\Omega}^r(\nabla e_0)
    +C_2\mathcal U_{\mathcal P}^1(R_{\mathrm H}(v)),    \\
    \widehat{\mathcal R}_{\mathrm W}^r
    (v;\mathcal P_\Omega,\mathcal P)
    &\coloneqq
    \alpha^{\mathrm W}\mathcal U_{\mathcal P_\Omega}^r(\nabla e_0)
    +\eta^{\mathrm W}\mathcal U_{\mathcal P_\Omega}^r(e_{\mathrm v})
    +\beta^{\mathrm W}\mathcal U_{\mathcal P}^1
        (R_{\mathrm W}(v)). 
\end{aligned}
\end{equation}
Then
\cref{cor:heat-residual-estimates,cor:wave-residual-estimate} and the
quadrature constructions above give
\begin{equation}\label{eq:residual-bound-chain}
    \mathcal N_\sigma(u_\sigma-v)
    \leq \mathcal R_\sigma(v)
    \leq \widehat{\mathcal R}_\sigma^r
    (v;\mathcal P_\Omega,\mathcal P), \qquad \sigma\in\{\mathrm H_0,\mathrm H_1,\mathrm W\},
\end{equation}
where $u_{\mathrm H_0}=u_{\mathrm H_1}=u_{\mathrm H}$.

First, we must ensure that our computable residual bounds satisfy the hypotheses of \cref{lem:verified-norm-bound-convergence} and hence converge to the true residual norms.

\begin{theorem}[Convergence of computable residual bounds]
\label{thm:computable-residual-bound-convergence}
Fix
$v=Bf$ as in \eqref{eq:boundary-ansatz}, and let
$\mathcal P_\Omega$ and $\mathcal P$ be as in \cref{assumption:rectangular-domain}.
Let
$\sigma\in\{\mathrm H_0,\mathrm H_1,\mathrm W\}$, $r\in\{0,1\}$, and suppose that the initial data for the PDE have the corresponding regularity specified in \cref{ass:data-derivative-bounds}.
Let $h_\Omega$ and $h_P$ be the maximal radii
of the cells in $\mathcal P_\Omega$ and $\mathcal P$, respectively.  There
is a constant $K>0$, independent of the partitions, such that
\begin{equation}\label{eq:computable-residual-bound-convergence}
    0\leq
    \widehat{\mathcal R}_{\sigma}^r
    (v;\mathcal P_\Omega,\mathcal P)-\mathcal R_{\sigma}(v)
    \leq K\left(h_\Omega^{r+1}+h_P^2\right).
\end{equation}
Consequently, along any sequence of such pairs of partitions for which
$h_\Omega,h_P\to0$,
\[
    \widehat{\mathcal R}_{\sigma}^r
    (v;\mathcal P_\Omega,\mathcal P)
    \longrightarrow \mathcal R_{\sigma}(v).
\]
\end{theorem}

\begin{proof}
Each $\widehat{\mathcal{R}}_\sigma^r$ is obtained by
replacing the norm of each residual $\phi$ in $\mathcal R_\sigma$ by its upper bound $\mathcal{U}^r(\phi)$ computed over the relevant spatial or spacetime partition. Thus it is sufficient to show that each $\mathcal{U}^r(\phi)$ converges at the required rate to the corresponding norm of $\phi$ as the partition is refined.  By \cref{lem:verified-norm-bound-convergence}, this follows if we can show that the computed $b_{q,C_k}(\phi)$ and $h_{q\ell,C_k}(\phi)$ for each residual are uniformly bounded (for the gradient residual $\nabla e_0$, the convergence then follows from the convergence of the bounds $\mathcal{U}^r(\partial_j e_0)$ for each component $\partial_j e_0$ and \cref{rem:vector-residual-bounds}).

The possible spatial and spacetime cell centres and radius vectors are contained in fixed compact sets since $\overline{\Omega_T}$ is compact. By \cref{ass:data-derivative-bounds}, the data bounds required by the initial data terms of $\widehat{\mathcal R}_\sigma^r$ are locally bounded and hence uniformly bounded. Since the boundary function $B$ is smooth and the moduli $\omega_B^\beta$ are continuous, the quantities $\mathsf B_k^\beta$ are similarly uniformly bounded. By \cref{lem:network-bound-local-boundedness}, the quantities $\mathsf F_k^\gamma$ are uniformly bounded, and hence the finite sums of products defining $\mathsf V_k^\alpha$ and $\mathsf V_{\mathrm v,k}^\alpha$ are as well. Thus, all the local bounds used in the initial data terms of $\widehat{\mathcal R}_\sigma^r$, as defined in \cref{eq:initial-displacement-coefficients,eq:initial-gradient-coefficients,eq:initial-velocity-coefficients}, are uniformly bounded. By \cref{lem:verified-norm-bound-convergence}, for each initial residual $\phi$ appearing in $\mathcal R_\sigma$, we have
\[
\begin{aligned}
\mathcal U_{\mathcal P_\Omega}^r(\phi)-\|\phi\|
&=O(h_\Omega^{r+1}).
\end{aligned}
\]

Finally, for the local bounds $h_{qr,C}^{\mathfrak p}$ defined for the PDE residuals in \cref{eq:appendix-pde-hessian-bound}, smoothness of $B$ and \cref{lem:network-bound-local-boundedness} give uniform boundedness since these bounds are sums of products of uniformly bounded quantities. By \cref{lem:verified-norm-bound-convergence}, this gives
\[
\mathcal U_{\mathcal P}^1(R_{\mathfrak p}(v))-\|R_{\mathfrak p}(v)\|=O(h_P^2),
\]
which completes the proof of \eqref{eq:computable-residual-bound-convergence}.
\end{proof}

\subsubsection{Local verification algorithm}

Given $v=Bf$, an estimate
$\sigma\in\{\mathrm H_0,\mathrm H_1,\mathrm W\}$, an initial data quadrature
order $r\in\{0,1\}$, and a desired accuracy $\eps_0>0$, we want to verify that
$v$ is within $\eps_0$ of the unknown solution $u_\sigma$ in the corresponding
error norm.  
Throughout, we assume the initial data regularity hypotheses of the residual estimate corresponding to $\sigma$ in \cref{cor:heat-residual-estimates,cor:wave-residual-estimate} and quadrature rule $Q^r$. 

\begin{enumerate}
    \item input: $v=Bf$, $\sigma\in\{\mathrm H_0,\mathrm H_1,\mathrm W\}$,
    $r\in\{0,1\}$, $\eps_0>0$, point values and moduli from
    \cref{ass:data-derivative-bounds,ass:B-moduli};
    \item initialise a spatial radius vector $\boldsymbol{\eps}_x\in(0,\infty)^d$
    and a time radius $\eps_t>0$ for which the spatial and spacetime
    partitions $\mathcal P_\Omega$ and $\mathcal P$ exist as in
    \cref{assumption:rectangular-domain};
    \item compute the initial data and PDE residual norm bounds over
    these partitions using the algorithms from \cref{sec:activation} and the
    constructions above, and combine them into
    $\widehat{\mathcal R}_{\sigma}^r
    (v;\mathcal P_\Omega,\mathcal P)$;
    \item if
    \begin{equation}\label{eq:halting-criterion}
        \widehat{\mathcal R}_{\sigma}^r
        (v;\mathcal P_\Omega,\mathcal P)<\eps_0,
    \end{equation}
    halt and return \textsc{yes}. Otherwise replace
    $(\boldsymbol{\eps}_x,\eps_t)$ by
    $(\boldsymbol{\eps}_x/2,\eps_t/2)$, refine the partitions accordingly,
    and start over from point~3.
\end{enumerate}

Even when the error of $v$ is smaller than $\eps_0$, there is in general no
guarantee that the local verification algorithm will ever halt, since the
residual estimate itself may be pessimistic.  However,
\cref{thm:computable-residual-bound-convergence} guarantees finite termination
whenever $\mathcal R_\sigma(v)<\eps_0$. This motivates the global algorithm.

\subsubsection{Global verification algorithm}

Fix $\sigma\in\{\mathrm H_0,\mathrm H_1,\mathrm W\}$ and $r\in\{0,1\}$.  Repeatedly using the local algorithm, we can formulate the
global algorithm that affirmatively answers \cref{question-sequence}. The explicit
construction above applies to the heat and wave equations, but the same
procedure works for any PDE admitting a residual framework of the form
\eqref{eq:XYZ-residuals} whose residual norms can be bounded by the
quadrature algorithms and those in \cref{sec:activation}.

\begin{enumerate}
    \item input: a tolerance $\eps_0>0$, the estimate $\sigma$, the initial data
    quadrature order $r$, and a sequence
    $\{v_m=Bf_m\}$ containing a subsequence $\{v_{m_k}\}$ converging to $u$
    in a space $Z$ for which
    \eqref{eq:XYZ-residuals} holds.
    \item set $n\gets1$;
    \item while no local algorithm instance has halted:
    \begin{enumerate}
        \item initialise the local algorithm for $v_n$;
        \item perform one iteration of the local algorithm for each
        $v_m$, $m=1,\dots,n$;
        \item set $n\gets n+1$.
    \end{enumerate}
    \item output the first index $m$ whose local run halts.
\end{enumerate}

The global algorithm halts after finitely many steps. Indeed, by \eqref{eq:XYZ-residuals}, there is an index $m_k$ for which
$\mathcal R_\sigma(v_{m_k})<\eps_0$.  For this fixed $v_{m_k}$,
\cref{thm:computable-residual-bound-convergence} implies that after finitely
many refinements the local algorithm satisfies
\eqref{eq:halting-criterion}. Thus it eventually halts and returns
an approximation $v_m$ satisfying $\norm{u-v_m}_X<\eps_0$.

\begin{remark}
The space $Z$ in \eqref{eq:XYZ-residuals} may be chosen according to the PDE
and the available regularity. For the heat and wave equations,
$Z=H^2(\Omega_T)$ is a natural candidate, since $\mathcal R_\sigma(v)
\le C \|u_\sigma-v\|_{H^2(\Omega_T)}$. In particular,
\citet[Theorem~5.1]{deryck2021approximation} shows that sufficiently regular
functions can be approximated by tanh networks in $W^{2,\infty}$. Thus if $q\coloneqq u/B$ is sufficiently regular, there are approximations $Bf_m\to u$ in $H^2$. However, approximation results of this kind do not guarantee that a given training procedure will compute such a sequence, and they cannot provide a computable error bound when $u$ is unknown. Hence, the  global algorithm above is required to verify a sufficiently accurate approximation in practice, if indeed it has been computed by the training procedure.
\end{remark}

\section{Numerical examples}\label{sec:numerical-examples}

\subsection{Heat and wave equation scenarios}

We consider the domain given by $T=1$ and $\Omega=(0,1)^d$. 
In both cases the initial displacement is
\[
    g(x)=\prod_{j=1}^d \sin(\pi x_j).
\]
For every spatial multi-index $\alpha$ with $|\alpha|\leq3$, the analytic form of $g$ provides the modulus
\[
    \omega_g^\alpha(r)
    =\min\{2\pi^{|\alpha|},d\pi^{|\alpha|+1}r\}.
\]
Consequently, the derivative bounds required by
\cref{ass:data-derivative-bounds} are
\[
 \mathcal D_g^\alpha(y,\boldsymbol\eps)
 \coloneqq
 |\partial^\alpha g(y)|+
 \omega_g^\alpha\p*{\norm{\boldsymbol\eps}_\infty},
 \qquad |\alpha|\leq3.
\]
We use initial velocity $h\equiv 0$, and hence
$\mathcal D_h^\alpha\equiv0$,
for the wave equation. The heat equation has $\kappa=0.1$ and the wave equation has $c=1$.
For this domain,
$\lambda_1=d\pi^2$ and
\[
C_{\Omega,d}
=\left(1+\frac{1}{d\pi^2}+\frac{1}{d^2\pi^4}\right)^{1/2}.
\]
Substituting $\kappa=0.1$ and $c=1$ into the constants in
\cref{eq:heat-alpha-beta-constants,eq:C-H1-H2,eq:wave-alpha-beta-constants}
gives
\[
\begin{gathered}
\alpha^{\mathrm H}=1+\sqrt{10}+\frac1{\sqrt{10}},
\qquad
\beta_d^{\mathrm H}=\frac{12+\sqrt{10}}{\pi\sqrt d},\\
C_{1,d}=1+\frac1{\sqrt{10}}+\sqrt{10}\,C_{\Omega,d},
\qquad
C_{2,d}=2+\sqrt{10}+10C_{\Omega,d},\\
\alpha^{\mathrm W}=\eta^{\mathrm W}=3,
\qquad
\beta_d^{\mathrm W}=3+\frac1{\pi\sqrt d}.
\end{gathered}
\]

\subsection{Quadrature results}\label{sec:separable-initial-displacement-data}
In each scenario, a model $v = Bf$ was trained with a standard PINN loss: discretised $L^2$ norms of the PDE and initial data residuals. No boundary loss is needed because of the factor $B$. The initial data residuals use a uniform spatial grid with radius $\eps_x=10^{-3}$, giving $500^d$ cells partitioning $(0,1)^d$. The PDE residuals use grids of size $500\times500$ for $d=1$ and $300^2\times750$ for $d=2$. For the most challenging case, $d=3$, the heat equation starts from a $100^3\times500$ base grid and uses adaptive refinement: each cell where the volume-normalised gap between the computed bound and midpoint value was larger than an empirically chosen threshold was uniformly subdivided into $2^{d+1}=16$ cells, before applying the same refinement procedure to the resulting cells once more.
For the three-dimensional wave equation, we used a fixed $200\times100^2\times500$ grid, i.e., with higher resolution in the first spatial dimension. These additional optimisations were found to improve quadrature accuracy for our trained models, but for different models or PDEs, other strategies may improve results.

Computed bounds for the initial displacement and initial displacement gradient residuals are given in \cref{tab:initial-displacement-residuals-q0-q1,tab:initial-displacement-gradient-residuals-q0-q1}, respectively. The affine rule $Q^1$ offers a significant improvement over the midpoint rule $Q^0$, but requires stronger assumptions on the data (which of course are satisfied in this case). For the initial velocity residuals in \cref{tab:wave-initial-velocity-residuals-q0-q1}, the improvement is marginal, mostly due to the fact that $h\equiv 0$. The final bounds in \cref{tab:energy-estimate-q1-master} are obtained by inserting the
$Q^1$ columns of the four preceding tables into
\cref{eq:heat-residual-estimate-l2,eq:heat-residual-estimate-h1,eq:wave-residual-estimate}, producing the cumulative bound for the LHS in each energy estimate.

Remarkably, our framework shows that even though a NN-based model is trained in a completely nonrigorous manner using only discretised $L^2$ losses for the PDE and initial data residuals, our verification algorithms provide rigorous error bounds in the substantially stronger norms displayed on the left-hand sides
of \cref{eq:heat-residual-estimates,eq:wave-residual-estimate}. We also point out that a natural application of our framework is to use the verification algorithms to guide training, e.g., by adding the residual bounds (or a proxy) to the loss function to encourage easier verification, or to focus training on more difficult regions of the domain. This is a direction of future work. The code and model weights used to produce the results are available at \url{https://github.com/emilhaugen/NetBounds}. AI coding assistance tools were used only to improve the readability and documentation of a pre-existing codebase\footnote{\url{https://github.com/emilhaugen/RigPINNs}} that was developed solely by one of the authors. All algorithms were designed and written by the authors.

\begin{table}[t]
\centering
\caption{Quadrature estimates and upper bounds for the initial displacement residual norm, as constructed in \cref{sec:initial-displacement-quad}.}
\label{tab:initial-displacement-residuals-q0-q1}
\small
\resizebox{\textwidth}{!}{%
\begin{tabular}{llrrr}
\toprule
Equation & $(d,L,w)$ & $\sqrt{Q_{\mathcal P_\Omega}^0(e_0)}$ & $\mathcal U_{\mathcal P_\Omega}^0(e_0)$ & $\mathcal U_{\mathcal P_\Omega}^1(e_0)$ \\
\midrule
Heat & $(1,2,128)$ & \num{7.0152e-05} & \num{3.0807e-03} & \num{7.3380e-05} \\
Heat & $(2,3,128)$ & \num{3.8112e-05} & \num{4.8498e-03} & \num{4.5751e-05} \\
Heat & $(3,4,128)$ & \num{4.9293e-05} & \num{6.3792e-03} & \num{6.4054e-05} \\
\midrule
Wave & $(1,2,256)$ & \num{4.2115e-04} & \num{3.3603e-03} & \num{4.2436e-04} \\
Wave & $(2,3,256)$ & \num{3.9864e-05} & \num{4.8515e-03} & \num{4.7537e-05} \\
Wave & $(3,3,256)$ & \num{1.1642e-04} & \num{6.4202e-03} & \num{1.2966e-04} \\
\bottomrule
\end{tabular}
}%
\end{table}

\begin{table}
\centering
\caption{Quadrature estimates and upper bounds for the initial displacement gradient residual norm, as constructed in \cref{sec:initial-displacement-gradient-quad}.}
\label{tab:initial-displacement-gradient-residuals-q0-q1}
\small
\resizebox{\textwidth}{!}{%
\begin{tabular}{llrrr}
\toprule
Equation & $(d,L,w)$ & $\left(\sum_j Q_{\mathcal P_\Omega}^0(\partial_j e_0)\right)^{1/2}$ & $\mathcal U_{\mathcal P_\Omega}^0(\nabla e_0)$ & $\mathcal U_{\mathcal P_\Omega}^1(\nabla e_0)$ \\
\midrule
Heat & $(1,2,128)$ & \num{7.2908e-04} & \num{1.2417e-02} & \num{7.3604e-04} \\
Heat & $(2,3,128)$ & \num{8.6705e-04} & \num{2.4283e-02} & \num{8.9859e-04} \\
Heat & $(3,4,128)$ & \num{1.3065e-03} & \num{3.7498e-02} & \num{1.3832e-03} \\
\midrule
Wave & $(1,2,256)$ & \num{1.6971e-03} & \num{1.3208e-02} & \num{1.7046e-03} \\
Wave & $(2,3,256)$ & \num{5.8680e-04} & \num{2.4118e-02} & \num{6.2037e-04} \\
Wave & $(3,3,256)$ & \num{1.6122e-03} & \num{3.7515e-02} & \num{1.6851e-03} \\
\bottomrule
\end{tabular}%
}
\end{table}

\begin{table}[t]
\centering
\caption{Quadrature estimates and upper bounds for the initial velocity residual norm, as constructed in \cref{sec:initial-velocity-quad}.}
\label{tab:wave-initial-velocity-residuals-q0-q1}
\small
\resizebox{\textwidth}{!}{%
\begin{tabular}{llrrr}
\toprule
Equation & $(d,L,w)$ & $\sqrt{Q_{\mathcal P_\Omega}^0(e_{\mathrm v})}$ & $\mathcal U_{\mathcal P_\Omega}^0(e_{\mathrm v})$ & $\mathcal U_{\mathcal P_\Omega}^1(e_{\mathrm v})$ \\
\midrule
Wave & $(1,2,256)$ & \num{4.7735e-04} & \num{4.8376e-04} & \num{4.7737e-04} \\
Wave & $(2,3,256)$ & \num{3.5101e-04} & \num{3.9677e-04} & \num{3.5126e-04} \\
Wave & $(3,3,256)$ & \num{2.2290e-04} & \num{2.9844e-04} & \num{2.2344e-04} \\
\bottomrule
\end{tabular}
}%
\end{table}

\begin{table}[t]
\centering
\caption{Quadrature estimates and upper bounds for the PDE residual norms, as constructed in \cref{sec:pde-residual-quad}.}
\label{tab:pde-residuals-q1-master}
\small
\resizebox{\textwidth}{!}{%
\begin{tabular}{llcrr}
\toprule
Equation & $(d,L,w)$ & Grid & $\sqrt{Q_{\mathcal P}^1(R_{\mathfrak p}(Bf))}$ & $\mathcal U_{\mathcal P}^1(R_{\mathfrak p}(Bf))$ \\
\midrule
Heat & $(1,2,128)$ & $500\times 500$ & \num{9.9469e-04} & \num{9.9673e-04} \\
Heat & $(2,3,128)$ & $300^2\times 750$ & \num{1.0057e-03} & \num{1.0212e-03} \\
Heat & $(3,4,128)$ & $100^3\times 500\,\mathrm{(+adaptive)}$ & \num{1.0021e-03} & \num{1.8986e-03} \\
\midrule
Wave & $(1,2,256)$ & $500\times 500$ & \num{6.3938e-03} & \num{6.4400e-03} \\
Wave & $(2,3,256)$ & $300^2\times 750$ & \num{6.6856e-03} & \num{7.0087e-03} \\
Wave & $(3,3,256)$ & $200\times 100^2\times 500$ & \num{3.0088e-03} & \num{1.0746e-02} \\
\bottomrule
\end{tabular}
}%
\end{table}

\begin{table}[t]
\centering
\caption{Energy-estimate bounds obtained by inserting the bounds from \cref{tab:initial-displacement-residuals-q0-q1,tab:initial-displacement-gradient-residuals-q0-q1,tab:wave-initial-velocity-residuals-q0-q1,tab:pde-residuals-q1-master} into \cref{eq:heat-residual-estimate-l2,eq:heat-residual-estimate-h1,eq:wave-residual-estimate}.}
\label{tab:energy-estimate-q1-master}
\small
\resizebox{\textwidth}{!}{%
\begin{tabular}{lllrrrr}
\toprule
Equation & $(d,L,w)$ & LHS & Initial displacement & Initial velocity & PDE residual & Verified LHS upper bound \\
\midrule
Heat ($L^2$ data) & \multirow{2}{*}{$(1,2,128)$} & $\mathcal N_{\mathrm H_0}(e_{\mathrm H})$ & \num{3.2864e-4} & -- & \num{4.8106e-3} & \num{5.1393e-3} \\
Heat ($H_0^1$ data) &  & $\mathcal N_{\mathrm H_1}(e_{\mathrm H})$ & \num{3.4228e-3} & -- & \num{1.5655e-2} & \num{1.9078e-2} \\
\addlinespace
Heat ($L^2$ data) & \multirow{2}{*}{$(2,3,128)$} & $\mathcal N_{\mathrm H_0}(e_{\mathrm H})$ & \num{2.0490e-4} & -- & \num{3.4851e-3} & \num{3.6900e-3} \\
Heat ($H_0^1$ data) &  & $\mathcal N_{\mathrm H_1}(e_{\mathrm H})$ & \num{4.0990e-3} & -- & \num{1.5753e-2} & \num{1.9852e-2} \\
\addlinespace
Heat ($L^2$ data) & \multirow{2}{*}{$(3,4,128)$} & $\mathcal N_{\mathrm H_0}(e_{\mathrm H})$ & \num{2.8687e-4} & -- & \num{5.2906e-3} & \num{5.5775e-3} \\
Heat ($H_0^1$ data) &  & $\mathcal N_{\mathrm H_1}(e_{\mathrm H})$ & \num{6.2703e-3} & -- & \num{2.9117e-2} & \num{3.5388e-2} \\
\midrule
Wave & $(1,2,256)$ & $\mathcal N_{\mathrm W}(e_{\mathrm W})$ & \num{5.1139e-3} & \num{1.4322e-3} & \num{2.1371e-2} & \num{2.7918e-2} \\
Wave & $(2,3,256)$ & $\mathcal N_{\mathrm W}(e_{\mathrm W})$ & \num{1.8611e-3} & \num{1.0538e-3} & \num{2.2604e-2} & \num{2.5519e-2} \\
Wave & $(3,3,256)$ & $\mathcal N_{\mathrm W}(e_{\mathrm W})$ & \num{5.0553e-3} & \num{6.7032e-4} & \num{3.4215e-2} & \num{3.9941e-2} \\
\bottomrule
\end{tabular}
}%
\end{table}

\FloatBarrier

\newpage
\appendix
\renewcommand{\theequation}{\thesection.\arabic{equation}}
\section{Boundary factor moduli and proofs for \cref{sec:activation}}\label{app:derivative-bound-proofs}

\subsection{Explicit moduli for the boundary factor}\label{app:boundary-factor-moduli}
For the boundary factor used in \cref{sec:applications}, 
\[
    B(x)=\prod_{i=1}^d x_i(1-x_i),
\]
write $s(r)=r(1-r)$, so that $B(x)=\prod_{i=1}^d s(x_i)$, and set
\[
    b_0\coloneqq \frac14,\qquad b_1\coloneqq 1,\qquad b_2\coloneqq 2,
    \qquad b_m\coloneqq 0\quad (m\geq 3).
\]
For any spatial multi-index $\beta\in\N_0^d$,
\[
    \partial^\beta B(x)
    =
    \prod_{i=1}^d s^{(\beta_i)}(x_i),
    \qquad
    \sup_{x\in(0,1)^d}|\partial^\beta B(x)|
    \leq
    \prod_{i=1}^d b_{\beta_i}.
\]
Thus $\partial^\beta B\equiv0$ whenever some $\beta_i\geq3$.  If $e_j$
denotes the $j$th coordinate multi-index, a valid nondecreasing modulus is
\begin{equation}\label{eq:explicit-B-moduli}
    \omega_B^\beta(r)
    \coloneqq
    r\sum_{j=1}^d\prod_{i=1}^d b_{\beta_i+\delta_{ij}},
    \qquad r\geq0.
\end{equation}
Terms with $\beta_i+\delta_{ij}\geq3$ vanish by the convention above.  Indeed,
for $x,y\in[0,1]^d$, the mean value theorem gives
\[
    |\partial^\beta B(x)-\partial^\beta B(y)|
    \leq
    \norm{x-y}_\infty
    \sum_{j=1}^d\sup_{z\in(0,1)^d}|\partial^{\beta+e_j}B(z)|,
\]
which is exactly \eqref{eq:explicit-B-moduli}.  In particular, this supplies
the moduli required by \cref{ass:B-moduli} for all $|\beta|\leq4$.

\subsection{Proofs delayed in \cref{sec:activation}}\label{app:xi-lemma-proof}
\begin{proof}[Proof of \cref{lem:xi-lemma}]
All operations are componentwise, so it suffices to prove the scalar statement.
Fix one component and write $\xi_i=\xi_{i,0}+\eta_i$ with $|\eta_i|\le E_i$.
Expanding the difference of products gives the exact identity
\[
    \prod_{i=1}^s(\xi_{i,0}+\eta_i)-\prod_{i=1}^s \xi_{i,0}
    =\sum_{\emptyset\ne I\subset\{1,\dots,s\}}
      \left(\prod_{i\in I}\eta_i\right)
      \left(\prod_{i\notin I}\xi_{i,0}\right).
\]
Taking absolute values and applying the triangle inequality gives
\[
\left|\prod_{i=1}^s \xi_i-\prod_{i=1}^s \xi_{i,0}\right|
\le
\sum_{\emptyset\ne I\subset\{1,\dots,s\}}
      \left(\prod_{i\in I}E_i\right)
      \left(\prod_{i\notin I}|\xi_{i,0}|\right).
\]
The right-hand side is precisely the expansion of
$\prod_i(|\xi_{i,0}|+E_i)-\prod_i|\xi_{i,0}|$.  This proves the scalar
inequality, and hence the componentwise vector inequality.
\end{proof}


\begin{proof}[Proof of \cref{lem:tanh-activation-bound}]\label{app:tanh-activation-bound-proof}
The statement is componentwise.  Fix one component and write
$t=z_j^k(y)$, $s=z_j^k(x)$, and $\rho=\rho_j^k$.  By assumption
$|s-t|\le\rho$.  Taylor's formula for $\mu^{(m)}$ at $t$ gives, for some point
$u$ between $s$ and $t$,
\begin{equation}\label{eq:tanh-taylor-expansion}
    \mu^{(m)}(s)-\mu^{(m)}(t)
    =\sum_{r=1}^{N_m}\frac{(s-t)^r}{r!}\mu^{(m+r)}(t)
    +\frac{(s-t)^{N_m+1}}{(N_m+1)!}\mu^{(m+N_m+1)}(u).
\end{equation}
We claim that for all $n\in \N$ and $v\in\R$,
\begin{equation}\label{eq:tanh-derivative-decay-bound}
    |\mu^{(n)}(v)|\le2^{n+1}n!e^{-2|v|}.
\end{equation}
For $n=1$, use the exponential representation of $\tanh$ and the relation $\mu' = 1 - \mu^2$ to see that
\[
    |\mu'(v)| = \frac{4}{(e^v+e^{-v})^2} = \frac{4e^{-2|v|}}{(1+e^{-2|v|})^2} \le 4e^{-2|v|}.
\]
Assuming the estimate holds through order $n$, the Leibniz rule and the fact that $|\mu(v)|\le1$ give
\[
    |\mu^{(n+1)}(v)|
    \le\sum_{k=0}^n\binom{n}{k}
        |\mu^{(k)}(v)\mu^{(n-k)}(v)|
    \le2^{n+2}(n+1)!e^{-2|v|}.
\]
This proves \eqref{eq:tanh-derivative-decay-bound} by induction.  Every point
between $s$ and $t$ satisfies $|u|\ge\max(0,|t|-\rho)$.  Substituting
\eqref{eq:tanh-derivative-decay-bound} into \eqref{eq:tanh-taylor-expansion} and using
$|s-t|\le\rho$ gives exactly the
corresponding component of
\eqref{eq:tanh-activation-bound}.
\end{proof}

\section{Explicit formulas for PDE residual bounds}
\label{app:initial-data-remainders}

Here we construct explicit cellwise bounds needed for the PDE residuals in
\cref{sec:pde-residual-quad}. We give the complete construction of $h_{qr,C}^{\mathfrak p}$ in
\eqref{eq:pde-residual-second-coefficients} for a spacetime cell $C$ with centre
$\zeta_C$ and radius vector $\boldsymbol{\eps}_C$. The construction uses a decomposition of the operator that retains the cancellations arising when the network $f$ nearly satisfies the PDE. Let $e_1,\dots,e_{d+1}$ be
the spacetime coordinate multi-indices, and set $\tau=e_{d+1}$. Write
$s_i=x_i(1-x_i)$, $\widehat B_i=\prod_{j\ne i}s_j$,
$H_i=(1-2x_i)\partial_i f-f$, and
$\lambda_{\mathrm H}=\kappa$, $\lambda_{\mathrm W}=c^2$.
Define
\begin{equation}\label{eq:appendix-pde-network-variation}
 \mathsf D_{f,C}^\alpha\coloneqq
 \begin{cases}
 \displaystyle\sum_{j=1}^{d+1}\eps_{C,j}
 \mathcal E_f^{e_j}(\zeta_C,\boldsymbol{\eps}_C),&|\alpha|=0,\\[1ex]
 \alg{|\alpha|}(\alpha,f,\zeta_C,\boldsymbol{\eps}_C),&|\alpha|\geq1.
 \end{cases}
\end{equation}
Then $$|\partial^\alpha f(\zeta)-\partial^\alpha f(\zeta_C)|
\leq\mathsf D_{f,C}^\alpha, \quad \zeta\in C.$$

For $j\leq d$, let
$I_{C,j}=[a_{C,j},b_{C,j}]
=[\zeta_{C,j}-\eps_{C,j},\zeta_{C,j}+\eps_{C,j}]\cap[0,1]$
and $\xi_{C,j}=\min\{b_{C,j},\max\{a_{C,j},1/2\}\}$. With
$s(x)=x(1-x)$, set
\begin{equation}\label{eq:appendix-pde-boundary-univariate-bounds}
 b_{C,j}^{(m)}\coloneqq\sup_{x\in I_{C,j}}|s^{(m)}(x)|
 =\begin{cases}
 s(\xi_{C,j}),&m=0,\\
 \max\{|1-2a_{C,j}|,|1-2b_{C,j}|\},&m=1,\\
 2,&m=2,\\
 0,&m\geq3.
 \end{cases}
\end{equation}
The boundary factor bounds, including time derivatives, are
\begin{equation}\label{eq:appendix-pde-boundary-product-bounds}
 \mathsf B_C^\beta\coloneqq
 \begin{cases}
 \displaystyle\prod_{j=1}^d b_{C,j}^{(\beta_j)},&\beta_{d+1}=0,\\
 0,&\beta_{d+1}>0,
 \end{cases}
 \qquad
 \widehat{\mathsf B}_{i,C}^\beta\coloneqq
 \begin{cases}
 \displaystyle\prod_{\substack{j=1\\j\ne i}}^d b_{C,j}^{(\beta_j)},
 &\beta_{d+1}=\beta_i=0,\\
 0,&\beta_{d+1}+\beta_i>0.
 \end{cases}
\end{equation}
Thus $\sup_C|\partial^\beta B|\leq\mathsf B_C^\beta$ and
$\sup_C|\partial^\beta\widehat B_i|\leq\widehat{\mathsf B}_{i,C}^\beta$.

For $\gamma\in\N_0^{d+1}$, define
\begin{equation}\label{eq:appendix-pde-operator-derivatives}
 \mathcal A_{\mathrm H}^\gamma
 \coloneqq\partial^{\gamma+\tau}f-\kappa\sum_{i=1}^d\partial^{\gamma+2e_i}f,
 \qquad
 \mathcal A_{\mathrm W}^\gamma
 \coloneqq\partial^{\gamma+2\tau}f-c^2\sum_{i=1}^d\partial^{\gamma+2e_i}f.
\end{equation}
The residual can be expanded as 
\begin{equation}\label{eq:pde-operator-aware-expansion}
 R_{\mathfrak p}(Bf)=B\mathcal A_{\mathfrak p}^0
 -2\lambda_{\mathfrak p}\sum_{i=1}^d\widehat B_iH_i,
 \qquad \mathfrak p\in\{\mathrm H,\mathrm W\}.
\end{equation}
The corresponding verified bounds are
\begin{equation}\label{eq:appendix-pde-operator-bounds}
 \begin{aligned}
 \mathsf A_{\mathrm H,C}^\gamma
 &\coloneqq|\mathcal A_{\mathrm H}^\gamma(\zeta_C)|
 +\mathsf D_{f,C}^{\gamma+\tau}
 +\kappa\sum_{i=1}^d\mathsf D_{f,C}^{\gamma+2e_i},\\
 \mathsf A_{\mathrm W,C}^\gamma
 &\coloneqq|\mathcal A_{\mathrm W}^\gamma(\zeta_C)|
 +\mathsf D_{f,C}^{\gamma+2\tau}
 +c^2\sum_{i=1}^d\mathsf D_{f,C}^{\gamma+2e_i}.
 \end{aligned}
\end{equation}
Moreover, set
\begin{equation}\label{eq:appendix-pde-H-derivatives}
 H_i^\gamma\coloneqq\partial^\gamma H_i
 =(1-2x_i)\partial^{\gamma+e_i}f-(2\gamma_i+1)\partial^\gamma f,
\end{equation}
and
\begin{equation}\label{eq:appendix-pde-H-bounds}
 \begin{aligned}
 \mathsf H_{i,C}^\gamma\coloneqq{}&|H_i^\gamma(\zeta_C)|
 +\bigl(|1-2\zeta_{C,i}|+2\eps_{C,i}\bigr)\mathsf D_{f,C}^{\gamma+e_i}
 +2\eps_{C,i}|\partial^{\gamma+e_i}f(\zeta_C)|\\
 &+(2\gamma_i+1)\mathsf D_{f,C}^\gamma
 \geq\sup_{\zeta\in C}|\partial^\gamma H_i(\zeta)|.
 \end{aligned}
\end{equation}

Finally, for $q,r=1,\dots,d+1$ and
$\mathfrak p\in\{\mathrm H,\mathrm W\}$, set
\begin{equation}\label{eq:appendix-pde-hessian-bound}
\begin{aligned}
 h_{qr,C}^{\mathfrak p}\coloneqq{}&
 \mathsf B_C^{e_q+e_r}\mathsf A_{\mathfrak p,C}^{0}
 +\mathsf B_C^{e_q}\mathsf A_{\mathfrak p,C}^{e_r}
 +\mathsf B_C^{e_r}\mathsf A_{\mathfrak p,C}^{e_q}
 +\mathsf B_C^{0}\mathsf A_{\mathfrak p,C}^{e_q+e_r}\\
 &+2\lambda_{\mathfrak p}\sum_{i=1}^d\Bigl(
 \widehat{\mathsf B}_{i,C}^{e_q+e_r}\mathsf H_{i,C}^{0}
 +\widehat{\mathsf B}_{i,C}^{e_q}\mathsf H_{i,C}^{e_r}
 +\widehat{\mathsf B}_{i,C}^{e_r}\mathsf H_{i,C}^{e_q}
 +\widehat{\mathsf B}_{i,C}^{0}\mathsf H_{i,C}^{e_q+e_r}\Bigr).
\end{aligned}
\end{equation}
The product rule and \eqref{eq:pde-operator-aware-expansion} give
$\sup_{\zeta\in C}|\partial_{qr}R_{\mathfrak p}(Bf)(\zeta)|\leq h_{qr,C}^{\mathfrak p}$.

\bibliographystyle{arxivpaper}
\bibliography{referencesClean}

\end{document}